\documentclass[final]{siamltex}
\usepackage{latexsym, graphicx, epsfig, amsmath, amsfonts}
\usepackage{verbatim}
\usepackage{subfigure}

\title{Equation-free, multiscale computation for unsteady random diffusion}
\author{Dongbin Xiu\thanks{Department of Chemical Engineering,
               Princeton University, Princeton, NJ 08544.}
             {(\tt{xiu@dam.brown.edu})}.
        \and Ioannis G. Kevrekidis
        \thanks{Department of Chemical Engineering, 
          Princeton University, Princeton, NJ 08544.}
        {(\tt{yannis@princeton.edu})}.}

\begin{document}

\maketitle

\begin{abstract}
We present an ``equation-free'' multiscale approach to the
simulation of unsteady diffusion in a random medium. 
The diffusivity of the medium is modeled as a random field with short
correlation length,
and the governing equations are cast in the form of stochastic 
differential equations.
A detailed fine-scale computation of such a problem requires 
discretization and solution of a large system of equations, 
and can be prohibitively time-consuming. 
To circumvent this
difficulty, we propose an equation-free approach, where
the fine-scale computation is conducted only for a (small) fraction of the
overall time. 
The evolution of a set of appropriately defined coarse-grained variables 
(observables)
is evaluated during the fine-scale computation, and   
``projective integration'' is used to accelerate the integration. 
The choice of these coarse variables is an important part of the approach:
they are the coefficients of pointwise polynomial expansions 
of the random solutions.
Such a choice of coarse variables allows us to reconstruct 
representative ensembles of fine-scale solutions with 
``correct" correlation structures, which is
a key to algorithm efficiency.
Numerical examples demonstrating accuracy and efficiency of
the approach are presented.
\end{abstract}
 

\pagestyle{myheadings}
\thispagestyle{plain}
\markboth{DONGBIN XIU AND IOANNIS KEVREKIDIS}
{EQUATION-FREE METHOD FOR RANDOM DIFFUSION}

\section{Introduction}

This paper is devoted to  numerical simulations of diffusion 
in a random medium whose material property, i.e., diffusivity (permeability,
conductivity), is characterized by small-scale, rough structures. 
This problem arises in the study of 
composite material properties, flow in porous media, etc
(e.g. \cite{Dagan89, Milton02}). 
Direct, fully resolved computations of the governing  equations in such media
can be prohibitively time consuming, 
as the fine-scale structures require discretizations resulting  to 
large degree-of-freedom calculations.
Hence, there has been a growing interest
in designing efficient alternative methods to solve the problems with 
the desired accuracy.

The properties of such media are typically modeled as 
{\em deterministic} smooth functions superimposed with fast oscillatory components.
One of the traditional approaches is to derive the {\it effective} properties 
of
such (heterogeneous) media -- the so-called ``homogenization" or ``upscaling"
techniques. 
These techniques are typically based on analytical asymptotic treatments
and  have been remarkably successful in several applications. 
Their applicability may, however, be restricted due to the assumptions 
that need to be made for the media
(See, for example, \cite{BensoussanLP78, CioranescuD99, E_CPA92}).
Numerical (as contrasted to analytical) approach to homogenization
are typically based on building multiscale basis functions into 
the spatial discretization. 
Methods along this
line of approach can be found in \cite{HouW_JCP97, HouWC_MC99,
HughesFMQ_CMAME98, MatacheBS_NM00, OberaiP_CMAME98} 
and are the subject of active research.

Alternatively, we can choose to model the medium properties as {\em random} 
fields,
to account for our insufficient knowledge and/or measurement error. 
For example, field data indicate that the
conductivity of many natural porous formations can be accurately
described by a lognormal distribution (e.g. \cite{Freeze_WRR75}).
The homogenization techniques developed for {\em deterministic media} can
be generalized to {\em random media}, and a comprehensive review can be found
in \cite{RenardM_AWR97}.

In addition to deriving equations describing an effective medium, as 
homogenization
does, many efforts have been devoted to direct detailed simulations 
of {\it random} media.
In this context, the corresponding governing equations,
e.g. the diffusion equation, Darcy's law, etc., are cast in the form
of stochastic equations and solved as such directly.
This approach further complicates the problem, since the governing equations
are now defined in (much) higher dimensional spaces, composed of both
the physical space and the space accounting for the parameterization of
the medium randomness (the random space).
The most straightforward numerical approach is the Monte Carlo method,
(e.g. \cite{Fishman96})
where repetitive deterministic simulations are conducted for particular 
realizations of the random functions describing the medium
properties, and what we are interested in here is the {\em statistics} of
the solution.
%
%
This approach, based on random sampling, 
can be computationally expensive because the
convergence rate of the ensemble averages, e.g., mean solution,
is relatively low (Monte Carlo simulations consisting
of $M$ realizations converge at a rate of $1/\sqrt{M}$.) 
There has been, therefore, a continuing interest in constructing non-sampling 
methods,
which include perturbation methods \cite{KleiberH92}, 
second-moment analysis \cite{LiuBM86}, stochastic Galerkin methods
\cite{XiuK_CMAME02, BabuskaTZ02, GhanemS91}, etc.
Among them, the stochastic Galerkin methods, also called ``generalized
polynomial chaos" expansions, 
have been  successful in many applications, when
the basis functions in the random space are appropriately chosen.
In particular, when the solution is  sufficiently smooth  
in the random space, stochastic Galerkin methods converge exponentially
fast \cite{XiuK_SISC02, XiuK_CMAME02, BabuskaTZ02}.

However, for the problem we intend to study in this paper -- random media
with short correlation length -- the stochastic Galerkin methods
become inefficient. This is because the short correlation length 
will induce a higher-dimensional random space, subsequently larger number
of equations to be solved -- the so-called ``curse-of-dimensionality''.
We remark that such a difficulty exists for all the existing
non-sampling methods.
On the other hand, although Monte Carlo methods are (formally) immune from
such a curse-of-dimensionality, their inherently slow convergence rate
can hardly improve the overall efficiency.

In this paper, we propose an equation-free, multiscale method to 
simulate diffusion in a random media with small-scale structure
(short correlation length).
The equation-free methods were first introduced in 
\cite{TheodoropoulosQK_PNAS00} and are designed to resolve multi-scale
problems efficiently. 
Such methods solve the equations for the effective, coarse-grained behavior 
without obtaining them in closed form; the quantities required for
these computations (residuals, action of Jacobians, time derivatives)
are {\it estimated} by solving the microscopic/stochastic model
with appropriately chosen initial conditions over short times
(and, in certain cases, only parts of the spatial domain).
The numerical results of these appropriately initialized short bursts
of microscopic computations are used 
to estimate the rate-of-change (or other quantities of interest)
of  appropriately defined {\it observables}: 
macroscopic variables characterizing the coarse-grained evolution.
These rates are then used by {\it projective integration} to evolve the
coarse-grained observables in time with (hopefully much)  larger time
steps (\cite{GearK_SISC03, GearK_JCP03, MartinezGK_JCP04}).
Thus, the time-consuming microscopic solvers are used for only
a (small) fraction of the overall time integration and 
no explicit knowledge of the equations governing the macroscopic
variables
is required (equation-free). 
This framework has been applied to a variety
of problems, ranging from bifurcation analysis of complex systems to
homogenization of periodic media
\cite{GearK_SISC03, GearLK_PLA03,KevrekidisGHKRT_CMS03,LiKGK_SIMMS03,
MakeevMK_JChP02, MakeevMPK_JChP02,RunborgTK_NL02, SamaeyRK_SIMMS04}.

For the random media with rough structure considered in this paper, we
employ the Monte Carlo method as the ``fine-scale'' solver. 
An orthogonal polynomial expansion of the fine-scale solution is conducted
(in principle) at every point in the physical space.
Our coarse-grained observables are the first few expansion coefficients of
such pointwise polynomial expansions on a relatively coarse grid;
the key assumption underlying our method is that {\it in principle} 
it should be possible to write a closed
equation for these observables that can successfully describe the
(long term) evolution of the solution statistics.
%

The particular assumption (observation) in this paper is
that, although
each individual realization of the solution is characterized by small spatial scales,
induced by the small scales in the diffusivity field, ensemble solution
averages  are characterized by larger, ``coarser"  scales.
That is, after possibly a short initial transient (relaxation), 
the ensemble solution
averages are significantly smoother than individual realizations in space,
and can therefore be accurately approximated with (hopefully significantly)
fewer degrees of freedom (e.g. on a coarse mesh).
This, in turn, implies that a closed equation exists (whether we can
explicitly derive it or not) for these averages on a coarser mesh;
this is precisely the equation that we will try to solve in this paper,
without explicitly deriving it.

The coefficients of the pointwise polynomial expansion of
the random solutions are representative of such ensemble averages,
and are observed to be smoother functions in space.
We thus expect that they can indeed be represented on a coarse mesh.
Since the explicit form of the 
governing equations for the evolution of such coefficients is unknown 
(to our best knowledge),
we employ the equation-free framework to compute with it.
In effect, we are trying to combine the simplicity of the Monte Carlo
implementation with the strengths of (generalized)
polynomial chaos representation:
instead of deriving and discretizing the equations for the appropriate
polynomial chaos coefficients via Galerkin expansion, we try to solve these 
equations through the design of ``just enough" short computational experiments
with the detailed direct solvers.
To this end, the rate-of-change of these coarse-grained variables are estimated
numerically from short bursts of fine-scale computation,
and propagated in time with larger steps via the projective integration technique.
The advantage of the present definition of the coarse variables is that
it allows us to reconstruct representative ensembles of fine-scale solutions 
with controlled accuracy. 
Numerical examples are presented to document the
accuracy and efficiency of the new algorithm.

The paper is organized as follows: in Section \ref{sec:formula} 
we formulate the
mathematical framework for diffusion in a random medium,
and subsequently the multiscale problem we will study. 
In Section 
\ref{sec:method}
we present the details of the construction of our ``equation-free''\
multiscale algorithm; in particular we focus on the ``lifting" step":
the construction of representative ensembles of fine scale solution
realizations.
Numerical examples are presented in Section 
\ref{sec:results},
and we conclude the paper with a general discussion in Section 
\ref{sec:summary}.

\section{Unsteady diffusion equations in random media}
\label{sec:formula}

In this section, we begin by presenting the mathematical framework for diffusion
in a random medium. 
We then formulate the multiscale problem that
we will study and discuss the difficulties in solving it efficiently.
 
\subsection{Formulations for random diffusion}

Let $D\in\mathbb{R}^d, d=1,2,3$ be a bounded polygonal domain,
$J=(0,T]\in \mathbb{R}^+=(0,\infty)$ for some fixed time $T>0$, and
$(\Omega,\mathcal{F},P)$ be a complete probability space. 
Here
$\Omega$ is the set of outcomes, $\mathcal{F}\subset 2^\Omega$ is the
$\sigma-$algebra of events and $P:\mathcal{F}\to[0,1]$ is a
probability measure. 
Let $D_T=D\times J$, and we study the 
following random diffusion equation: 
find a stochastic function,
$u:\Omega\times\bar{D}_T\to\mathbb{R}$, 
such that for $P-$almost everywhere $\omega\in\Omega$, 
the following equation holds:
\begin{equation} \label{eqn}
\left\{
\begin{array}{rcll}
u_t(\omega,\cdot) - 
\nabla\cdot\left(\kappa(\omega,x)\nabla u(\omega,\cdot)\right)
&=&f(\omega,\cdot), & \textrm{in } D_T, \\
u(\omega,\cdot)&=&0, & \textrm{on } \partial D\times [0,T],\\
u(\omega,\cdot)&=& u_0(\omega,x), & \textrm{on } D\times\{t=0\},
\end{array}
\right.
\end{equation}
where $u$ is the unknown and 
$u_t=\partial u/\partial t$ its time derivative.  
$\kappa, u_0:\Omega\times D\to\mathbb{R}$, and
$f:\Omega\times D_T\to\mathbb{R}$ are the known stochastic functions
with continuous and bounded covariance functions. 
Denote by $B(A)$ the
Borel $\sigma-$algebra generated by the open subsets of $A$, then
$\kappa$ and $u_0$ are assumed to be measurable with the $\sigma-$algebra 
$(\mathcal{F}\otimes B(D))$ and $f$ with
$(\mathcal{F}\otimes B(D_T))$. 

The following assumptions are made on the input stochastic data:
\begin{enumerate}
\item $\kappa$ is bounded and uniformly coercive, i.e.
     \begin{equation}  \label{coercivity}
      \exists~ \kappa_{\textrm{min}}, 
      \kappa_{\textrm{max}}\in (0,+\infty):
      P\left(\omega\in\Omega: \kappa(\omega,x)\in [\kappa_{\textrm{min}},
      \kappa_{\textrm{max}}], \forall x\in\bar{D}\right)=1.
     \end{equation}
Also, $\kappa$ has a uniformly bounded and continuous first
derivative, i.e. there exists a real deterministic constant $C$ such
that
\begin{equation} \label{k_L1}
P\left(\omega\in\Omega: \kappa(\omega,\cdot)\in C^1(\bar{D})
\textrm{ and }
\textrm{max}_{\bar{D}}\left|\nabla_x
  \kappa(\omega,\cdot)\right|<C\right)=1.
\end{equation}

\item $f\in L^2(\Omega)\otimes L^2(D_T)$, i.e.
\begin{equation}
\int_\Omega\int_J\int_D f^2(\omega,x,t)dxdt dP(\omega) < +\infty.
\end{equation}

\item $u_0\in L^2(\Omega)\otimes L^2(D)$, i.e.
\begin{equation}
\int_\Omega\int_D u_0^2(\omega,x)dx dP(\omega) < +\infty.
\end{equation}

\end{enumerate}

\subsection{Finite-dimensional noise and variational form}
\label{sec:weak}

For the problem \eqref{eqn} to be practically solvable numerically, 
it should be possible to reduce the
infinite-dimensional probability space to a finite-dimensional
space. 
This can be accomplished by characterizing the probability
space by a finite number of random variables.
Such a procedure,  termed as
the ``finite-dimensional noise assumption'' in \cite{BabuskaTZ02},
is often achieved
via a certain type of decomposition which can approximate the
target random process with desired accuracy. 
One of the choices is 
the Karhunen-Lo\'eve type expansion \cite{Loeve77},
which is based on the spectral decomposition of
the covariance function of the input random process.
(e.g. \cite{GhanemS91, XiuK_CMAME02}). 
Following a decomposition and assuming that the random inputs can
be characterized by $N$ random variables, we can rewrite
the random inputs in the following abstract form,
\begin{equation} \label{fd_noise}
\kappa(\omega,x)=\kappa(Y_1(\omega),\dots,Y_N(\omega),x),
\textrm{ and }
f(\omega,x,t)=f(Y_1(\omega),\dots,Y_N(\omega),x,t),
\end{equation}
where $N\geq 1$ is a finite integer and
$\{Y_n\}_{n=1}^N$ are real random variables with zero mean
value, unit variance, and their images $\Gamma_{n,N}\equiv
Y_n(\Omega)$ are bounded intervals in $\mathbb{R}$ for $n=1,\dots,N$.
Moreover, we assume that each $Y_n$ has a density function
$\rho_n:\Gamma_{n,N}\to\mathbb{R}^+$ for $n=1,\dots,N$, and denote
$\rho(y),\forall y\in\Gamma$  the joint probability density of
$(Y_1,\dots,Y_n)$ and
$\Gamma\equiv\prod_{n=1}^N\Gamma_{n,N}\subset\mathbb{R}^N$
the support of such density. 
The expectation operator is subsequently defined as $\mathbb{E}(f)=
\int_\Gamma f(y)\rho(y)dy$.
Note when random variables
$Y_n,n=1\dots,N$ are independent, we have
$\rho(y)=\prod_{n=1}^N\rho_n(y_n),\forall y\in\Gamma$. 

After the finite-dimensional characterization of the random inputs
\eqref{fd_noise}, the unsteady
diffusion equation (\ref{eqn}) can be expressed in the following
strong form,
\begin{equation} \label{fd_strong}
\left\{
\begin{array}{rcll}
u_t(y,x,t)-\nabla\cdot\left(\kappa(y,x)\nabla u(y,x,t)\right)
&=&f(y,x,t),  &\forall (y,x,t)\in\Gamma\times D\times J,\\
u(y,x,t)&=& 0, &\forall (y,x,t)\in\Gamma\times\partial D\times [0,T],\\
u(y,x,0)&=&u_0(y,x), &\forall (y,x)\in\Gamma\times D.
\end{array}
\right.
\end{equation}

Often we seek its weak solution satisfying the
following variational form:
find $u\in L_\rho^2(\Gamma)\otimes L^2(0,T; H_0^1(D))$ with
$u_t\in  L_\rho^2(\Gamma)\otimes L^2(0,T; H^{-1}(D))$ such that
\begin{equation} \label{fd_weak}
\left\{
\begin{array}{ll}
 \mathcal{I}_\rho(u_t,v)
+\mathcal{K}_\rho(u,v;\kappa)
= \mathcal{I}_\rho(f,v),
& \forall v\in L_\rho^2(\Gamma)\otimes H_0^1(D),\\
u(t=0)=u_0
& ~
\end{array}
\right.
\end{equation}
where
\begin{displaymath}
\mathcal{I}_\rho(v,w) = \int_\Gamma\rho(y)\int_D v(y,x)w(y,x) dxdy,
\end{displaymath}
and
\begin{displaymath}
\mathcal{K}_\rho(v,w;\kappa)=
\int_\Gamma\rho(y)\int_D\kappa(y,x)\nabla v(y,x)\cdot\nabla w(y,x)dxdy.
\end{displaymath}

Note that problem (\ref{fd_strong}) or (\ref{fd_weak}) becomes 
a $(N+d)$ dimensional
problem, where $d$ is the dimensionality of the physical space
$D$ and $N$ is the
dimensionality of the random space $\Gamma$.

\subsection{Formulations for multiscale problems}

In this section we formulate 
the multiscale problem associated with the stochastic
diffusion equation (\ref{eqn}). 
In particular, we consider the
problem where the random inputs, e.g. diffusivity $\kappa$, have
very small correlation length $l_\kappa\ll 1$, 
compared to the (macroscopic) domain of interest $D_T$.
For notational convenience, hereafter we restrict our exposition to
problems with only $\kappa$ being the random input, and study,
for $P$-almost
everywhere $\omega\in\Omega$,
\begin{equation} \label{multi_eqn}
\frac{\partial u^\varepsilon}{\partial t}(\omega,\cdot)=
\nabla\cdot\left[\kappa\left(\omega,\frac{x}{\varepsilon}\right)\nabla
u^\varepsilon(\omega,\cdot)\right] + f(x),\qquad \textrm{in } D_T,
\end{equation}
\begin{equation} \label{multi_bc}
u^\varepsilon(\omega,\cdot)=0, \qquad \textrm{on }\partial D\times[0,T],
\end{equation}
\begin{equation} \label{multi_ic}
u^\varepsilon(\omega,\cdot)=u_0^\varepsilon(x) 
\qquad \textrm{on }D\times\{t=0\}.
\end{equation}
Here we have assumed that the diffusivity $\kappa$ satisfies the conditions
(\ref{coercivity}) and (\ref{k_L1}), and is a homogeneous random 
field with a short correlation length 
$l_\kappa\sim O(\varepsilon)\ll O(1)$, i.e.
\begin{equation} \label{corr}
C_{\kappa}(x,y) = C(|x-y|/ l_\kappa),\qquad l_\kappa\sim
O(\varepsilon)\ll O(1),
~(x,y)\in \bar{D}\times \bar{D},
\end{equation}
where 
$C_{v}(x,y)\equiv\mathbb{E}\left[
(v(\omega,x)-\mathbb{E}[v(\omega,x)])
(v(\omega,y)-\mathbb{E}[v(\omega,y)])\right]$
is the two-point covariance function and $l_\kappa$ is the correlation length.
Again we characterize the diffusivity field
$\kappa(\omega,x/\varepsilon)$  by $N$ independent
random variables as in (\ref{fd_noise}). Hence, problem
(\ref{multi_eqn}) is in $(N+d)$ dimensions, and we can formulate it
in the strong and weak forms similar to 
(\ref{fd_strong}) and (\ref{fd_weak}), respectively.
 
The discretization in the
spatial domain $D\subset\mathbb{R}^d$ can be conducted via any standard
technique, e.g., finite difference, finite elements, etc., with a
maximum mesh spacing parameter $\delta>0$. To fully resolve
(\ref{multi_eqn}), we need to employ a fine mesh with $\delta<\varepsilon$. 
From a numerical point of view, very small mesh spacing $\delta$ often
results in restrictions on the size of time steps of numerical
schemes, and such restrictions are particularly severe for explicit schemes.
Hence, a fine-scale computation of (\ref{multi_eqn}) requires
computations with very small time steps on a very fine mesh. 

The discretizations in the $N-$dimensional random space $\Gamma$
can be conducted in different ways. The recently developed
stochastic Galerkin methods, or generalized polynomial chaos,
are extensions of the classical polynomial chaos which is based on the
Wiener-Hermite expansion \cite{GhanemS91}. These extensions include
the non-Hermite global orthogonal polynomial expansion from the Askey
scheme \cite{XiuK_SISC02, XiuK_CMAME02, XiuK_IJHMT03}, piecewise
polynomial expansions \cite{BabuskaC_CMAME02, BabuskaTZ02, DebBO01},
and wavelet basis \cite{MaitreKNG_JCP04,MaitreNGK_JCP04}.
The stochastic Galerkin methods have fast convergence as the polynomial
order is increased. In fact, under sufficient regularity requirements,
exponential convergence has been proved for stochastic elliptic equations
in \cite{BabuskaTZ02} and shown numerically for various stochastic 
equations in \cite{XiuK_SISC02, XiuK_CMAME02, XiuK_JCP03}. 
The total number of expansion terms, $K$, however, depends not only on
the order of polynomials but also the dimensionality $N$ of the random space.
When $N\gg 1$ is very large, $K$ 
increases fast with increasing order of polynomials. This significantly
reduces the convergence rate with respect to the number of expansion
terms for stochastic Galerkin methods. 
To this end, it may be necessary to resort to the Monte Carlo
method, as its convergence rate, $1/\sqrt{M}$ where $M$ is the number
of realizations, albeit slower, is independent of the value of $N$.

The number of random variables, $N$, used to represent the random process
$\kappa(\omega,x)$ depends on, among other factors, the correlation
length of $\kappa$. 
Although one may choose different decomposition
methods, in general, the value of $N$ is inversely proportional to
the correlation length $l_\kappa$.
For the problem we consider here, $\l_\kappa\ll O(1)$ implies
$N\gg 1$, and problem \eqref{multi_eqn} is in a high-dimensional
random space. Subsequently,
the fine-scale computation of (\ref{multi_eqn})
requires 
a large number of discretization terms in the $N$-dimensional 
random space $\Gamma$ (by either a large number of expansion terms
from a stochastic Galerkin method at a given order, 
or a large number of realizations
from a Monte Carlo method),
a fine mesh in the physical space $D$ to resolve the small scales induced
by $l_\kappa\ll O(1)$, and very small time steps in the time domain $J$.
Such computations can be prohibitively time-consuming.

\section{An equation-free multiscale method}
\label{sec:method}

In this section we present an equation-free multiscale method for
the integration of equation 
(\ref{multi_eqn}); other tasks (such as fixed point
algorithms for its stationary states) can also be formulated in 
an equation-free context (see the discussion in Section \ref{sec:summary}).
The key feature of the method is that the costly
fine-scale computations of \eqref{multi_eqn} 
are only conducted for a small fraction of the total integration time.
During the fine-scale computations, the rate of change of 
appropriate coarse-scale
variables is estimated numerically and subsequently represented
on a coarser mesh and integrated in
time with large time steps. 
The choice of such coarse variables is based upon the following
assumption (observation):
although each individual realization of the solution of (\ref{multi_eqn})
is characterized by small spatial scales, the ensemble averages, e.g. moments, 
of the solutions are significantly smoother in space (characterized
by much larger scales), i.e.
\begin{equation} \label{SSS}
\mathbb{E}\left[g\left(
u^\varepsilon(\omega,\frac{x}{\varepsilon},t)\right)
\right] = U_g(x,t), \qquad \forall g\in \mathcal{C},
\end{equation}
where $\mathcal{C}$ denotes a set of smooth functions.
Variables $U_g(x,t)$ are the ``coarse-grained''
variables, and we will describe
in detail their construction in the following Section.
Hereafter, we will drop the superscript $\varepsilon$ of the
fine-scale variables $u^\varepsilon$.

Equation (\ref{multi_eqn}) defines an evolutionary process
\begin{equation} \label{micro_eqn}
\frac{\partial u}{\partial t}(\omega,x,t)= 
r(u), \qquad (\omega,x)\in\Omega\times D,
\end{equation} 
characterized by a solution operator
$\{s(t)\}$, which forms a semigroup $u(\cdot,t)=s(t)u(\cdot,0)$ in
$t\in\ J$. 
We will assume that the set of properly defined coarse-scale variables from
(\ref{SSS}) satisfy, possibly after a short transient (relaxation) period,
{\it closed} differential equations
\begin{equation} \label{macro_eqn}
\frac{\partial U_g}{\partial t}(x,t) = R(U_g), \qquad (x,t)\in D_T.
\end{equation}
Note that typically there are multiple number of coarse variables. Subsequently,
$U_g$ is a vector field and \eqref{macro_eqn} is a system of equations.
We also remark that the explicit knowledge of the equation
(\ref{macro_eqn}) may not be easy to obtain, or may be too complicated
to be of any practical use if it were known.

\begin{figure}[htbp]
   \centerline{
   \psfig{file=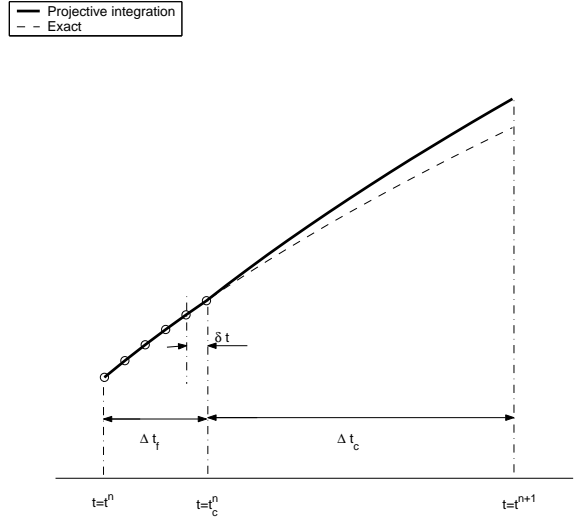,width=3in}}
\caption{Sketch of the multiscale equation-free integration over one
         global time step.}
\label{fig:graph}
\end{figure}

The general procedure for the equation-free projective integration methods 
(see, for example, \cite{GearK_SISC03, KevrekidisGHKRT_CMS03, ShortManifesto}) 
over one global time step $\Delta t$, starting
at $t=t^n$ and ending at $t=t^{n+1}$, consists of the following
key components:
\begin{itemize}
\item a ``restriction'' operator $\mathcal{P}$ 
      to evaluate the coarse-grained variables
      from the ensemble of fine-scale computations, i.e. $U_g=\mathcal{P}u$,
      and a ``lifting'' operator $\mathcal{Q}$
      to construct the representative ensemble of
      fine-scale solutions from the coarse-scale variables, i.e.
      $u=\mathcal{Q}U_g$;
\item $n_f>1$ steps of fine-scale computations of (\ref{micro_eqn}) 
with a small time step $\delta t$, 
where we will define $\Delta t_f=n_f\delta t$ and an intermediate
time level $t^n_c=t^n+\Delta t_f$;
\item one step of coarse projective integration of the coarse-scale
equation (\ref{macro_eqn}) with a time step of the size
$\Delta t_c=n_c\delta t, n_c>1$.
\end{itemize}
Since $\Delta t_c$
is associated with the time scale of the coarse variables $U_g$
defined in (\ref{SSS}),
we usually have $n_c\gg 1$. The global time step is $\Delta t=t^{n+1}
-t^n=\Delta t_f+\Delta t_c=(n_f+n_c)\delta t$. Figure \ref{fig:graph}
is a graphical illustration of the notations.

Specifically, for the multiscale diffusion problem in a random medium
described by (\ref{multi_eqn}), the 
equation-free integration over one global time step $\Delta t$
consists of the following steps: 
Given a fine and coarse computational mesh,
and $U_g^n(x)\equiv U_g(x,t^n)$ on the {\em coarse mesh}, 
\begin{enumerate}
\item {\em Lifting} (or {\em Reconstruction}): 
      Generate
      an ensemble of random solutions $u^{n}(\omega,x)\equiv
      u(\omega,x,t^n)$ on the fine
      mesh. 

\item {\em Fine-scale computation}: Fully resolve 
      (\ref{multi_eqn}) by using the fine mesh in $D$, a small time
      step $\delta t$ in $J$, and an appropriate
      method in $\Gamma$ (e.g., the Monte Carlo  method). 
      Such a fine-scale integration is conducted
      only for a short period of
      time, from $t^n$ to the intermediate time $t^n_c$, 
      i.e. $u(\cdot,t)=s(t)u(t^n)$, for $t^n\leq t\leq t^n_c=
      t^n+n_f\delta t$. Here $n_f\geq 1$ such that $n_f\delta
      t_f\sim t_R\ll t_M$, where $t_R$ is the local relaxation time of
      the fine-scale process and it is assumed to be much shorter
      than $t_M$, the typical time scales of the coarse-scale process
       (\ref{macro_eqn}).

\item {\em Restriction}: Evaluate the coarse variables $U_g(t)$ defined in
  (\ref{SSS}) on the coarse mesh, for $t^n\leq t\leq t^n_c$. 

\item {\em Coarse-scale integration}: Estimate the time derivatives of the
      coarse variables $U_g$ at $t=t^n_c$ and integrate 
      the coarse-scale equations
  (\ref{macro_eqn}) to $t^{n+1}$ via the projective integration method 
  {\em on the coarse mesh}, 
  with a time step $\Delta t_c=n_c\delta t=t^{n+1}-t^n_c$.
  Here $\Delta t_c\sim t_M\gg t_R$.

\end{enumerate}

We now present a detailed description of each of the steps, starting
with the more straightforward step -- the fine-scale computation.

\subsection{Fine-scale computation} \label{sec:mcfem}

The objective of the fine-scale computation is to fully resolve 
(\ref{multi_eqn}). To this end, any conventional spatial and
temporal discretization scheme can be employed, e.g. finite
difference or finite elements. Since the
dimensionality is high ($(N+d)$-dimensional), we employ a 
Monte Carlo simulation (MCS) in the random space $\Gamma$. Here we
illustrate the formulation of a Monte Carlo finite element method
(MCFEM). 

Denote $X_\delta^d\subset H_0^1(D), D\subset\mathbb{R}^d$ a family of
piecewise linear finite element approximation spaces, with a maximum
mesh spacing parameter $\delta>0$. This is the {\em fine mesh}, as
we choose $\delta<O(\varepsilon)\ll O(1)$ to fully resolve all spatial
scales.
We shall assume all the standard assumptions
on the finite element triangulation, and its approximation estimate,
i.e.
\begin{equation} \label{fem_approx}
\min_{\chi\in X_\delta^d}\Arrowvert v-\chi\Arrowvert \leq C \delta
\Arrowvert v\Arrowvert_{H^2(D)}, \qquad
\forall v\in H^2(D)\cap H_0^1(D),
\end{equation}
where $C>0$ is a constant independent of $v$ and $\delta$.

The MCFEM formulation for (\ref{multi_eqn}) is as following:
\begin{itemize}
\item Prescribe the number of realizations, $M$, and a piecewise linear finite
element space on $D$, $X_\delta^d$, defined as above.
\item For each $j=1,\dots,M$, sample i.i.d. realizations of the 
diffusivity $\kappa(\omega_j,\cdot)$ and find the corresponding
approximation $u_\delta(\omega_j,\cdot)\in L^2(0,T;X_\delta^d)$ 
with $\frac{\partial u_\delta}{\partial t}\in L^2(0,T;(X_\delta^d)')$
such that
\begin{equation} \label{mcfem}
\left(\frac{\partial u_\delta}
{\partial t}(\omega_j,\cdot), \chi\right)_\delta+
\int_D \kappa(\omega_j,\cdot)\nabla u_\delta(\omega_j,\cdot)\cdot
\nabla \chi dx=(f,\chi)_\delta, \quad \forall \chi\in X_\delta^d, t\in J,
\end{equation}
where $(\cdot,\cdot)_\delta$ is the usual inner product in $X_\delta^d$.
\item Process the solution ensemble to generate statistics, e.g.
      $\mathbb{E}(u) = 
       \frac{1}{M}\sum_{j=1}^M u_\delta(\omega_j,\cdot)$.
\end{itemize}

For more detailed discussion on the stochastic finite element spaces for
elliptic problems, see \cite{BabuskaTZ02}; for numerical examples and
implementations of stochastic Galerkin methods for steady/unsteady diffusion
equations, 
see \cite{XiuK_CMAME02, XiuK_IJHMT03}.

\subsection{Restriction and Lifting} \label{sec:projection}

The {\em restriction} from the fine-scale variables $u$ 
to coarse-scale variables $U_g$ 
consists of two steps: ``random restriction'' and ``spatial restriction''.
First, the fine-scale variables
$u$ are averaged to $U_g$ in the random space according to
(\ref{SSS}) (random restriction);
then the coarse variables $U_g$ are further restricted from the fine
mesh to the coarse mesh (spatial restriction) justified by the assumption/observation that they
are smoother. 
The {\em lifting procedure}  is the reverse of the restriction.
We now describe
the details of the two procedures in both the random space
$\Gamma$ and the physical space $D$.

\subsubsection{Operations in random space} \label{sec:I_rand}

The solution of the fine-scale computation via MCFEM
in section \ref{sec:mcfem}, or other effective methods, generates an
ensemble of $M$ realizations of the
random solution $u_\delta$ at any $x\in X_\delta^d$
and $t\in J$. 
For any fixed $(x,t)$, $u_\delta(\omega,\cdot)$ is a random variable, and
we seek to  represent such a random variable by an orthogonal
polynomial approximation and define the expansion coefficients as our
coarse-grained observables (variables). 
Hereafter we drop the subscript $\delta$ in
$u_\delta$, the fine-scale numerical solution of $u$,
and seek to approximate it by $\mathcal{I}_\omega u$ for any fixed $(x,t)$,
\begin{equation} \label{proj_oper}
u(\omega,x,t)\simeq\mathcal{I}_\omega u(\omega,x,t)
=\sum_{k=0}^{K} \bar{U}_k(x,t)\Phi_k(\xi(\omega)),
\end{equation}
where $\{\Phi_k(\xi(\omega))\}_{k=0}^{K}$ is a set of orthogonal polynomials
$\{\Phi\}$ in term of random variable $\xi(\omega)$. 
The expansion coefficients are determined by
\begin{equation} \label{proj_coef}
\begin{array}{rcl}
\bar{U}_k &=& \frac{1}{\int_\Omega\Phi_k^2 dP}
\int_\Omega u(\omega,\cdot)\Phi_k(\xi(\omega))dP(\omega)\\
 &=&
\frac{1}{\mathbb{E}\left[\Phi_k^2\right]}
\mathbb{E}\left[u(\omega,\cdot)\Phi_k(\xi(\omega))
\right], \quad k=0,\cdots,K,
\end{array}
\end{equation}
where the orthogonality of the basis functions has been used. 
Such a representation
of random variables is commonly used in practice
(cf. \cite{GhanemS91, XiuK_SISC02, XiuK_CMAME02}). The 
correspondence between the type of orthogonal polynomials $\{\Phi\}$ 
and the type of random
variable $\xi(\omega)$ includes Hermite-Gaussian, Jacobi-beta, 
Laguerre-gamma, etc (see \cite{XiuK_SISC02} for details). 
The convergence of such orthogonal polynomial expansions is assumed to
be of the form
\begin{equation} \label{ms_conv}
\Arrowvert u(\omega,\cdot) - \mathcal{I}_\omega u(\omega,\cdot)\Arrowvert
_{L^2(\Omega)}
=O(K^{-\gamma_\omega})\to 0, \qquad K\to\infty, ~\gamma_\omega>0,
\end{equation}
where we have assumed the convergence rate scales as $K^{-\gamma_\omega}$ for
some positive number $\gamma_\omega>0$, which depends on the
smoothness of $u(\omega)$.
We remark that a complete theoretical analysis on the convergence of
different basis remains an open issue.
For numerical examples of the approximations of a random variable via
different sets of basis, see \cite{XiuK_SISC02}.

The expansion coefficients $\{\bar{U}\}_{k=1}^K$ are the ensemble
averages of $u$, as defined in (\ref{proj_coef}), and under assumption
(\ref{SSS}), they become the coarse variables with larger
spatial scale, i.e., with smoother profiles in the physical space $D$. 
For instance, the
coefficient $\bar{U}_0(x,t)$ is the mean field of $u$ and is in general
smooth.

The finite-term polynomial approximation \eqref{proj_oper}
defines two operations between
$u$ and $\{\bar{U}_k\}$, i.e., the ``restriction'' operator in
random space $\mathcal{P}_\omega$ such that,
\begin{equation} \label{P_oper}
\{\bar{U}_k(\cdot)\}_{k=0}^K = \mathcal{P}_\omega u(\omega,\cdot) ,
\end{equation}
and the ``lifting'' operator 
$\mathcal{Q}_\omega$ such that
\begin{equation} \label{Q_oper}
\mathcal{I}_\omega u(\omega,\cdot)=\mathcal{Q}_\omega
\{\bar{U}_k(\cdot)\}_{k=0}^K,
\end{equation}
where operation $\mathcal{P}_\omega$ is accomplished by (\ref{proj_coef}), 
and $\mathcal{Q}_\omega$  by generating random samples of $\xi(\omega)$ in
equation (\ref{proj_oper}). Obviously both $\mathcal{P}_\omega$ 
and $\mathcal{Q}_\omega$
are linear operators, $\mathcal{Q_\omega P_\omega}=\mathcal{I}_\omega$ and 
$\mathcal{I}_\omega\mathcal{I}_\omega=\mathcal{I}_\omega$.

We remark that the expansion (\ref{proj_oper}) is different from the
traditional polynomial chaos expansion. 
Expansion (\ref{proj_oper})
is  a pointwise approximation at fixed locations in $(x,t)$,
and hence only requires a one-dimensional (in the random space) polynomial
basis of $\{\Phi_k(\xi(\omega))\}$ where the random variable
$\xi(\omega)$ associated with the basis is different at different
locations in the physical space. 
On the other hand, the traditional polynomial chaos
expansion is written in the full $N-$dimensional random space for {\em all}
locations of $(x,t)$. While operators (\ref{P_oper}) and
(\ref{Q_oper}) are one-dimensional in the random space, 
they do not offer us an easy way
to obtain the governing equations for the coarse variables $\bar{U}$,
shown as in (\ref{macro_eqn}). (On the other hand, we can readily
derive the governing equations in the $N-$dimensional random space via a Galerkin
method if the traditional polynomial chaos expansion is employed.)
We will show in the next section that we can circumvent the difficulty
of not having the governing equations by using the ``equation-free''
approach.

\subsubsection{Operations in physical space} \label{sec:I_phys}

Since we have assumed the coarse variables $\{\bar{U}_k(x,t)\}$ 
are smooth in space, they can be accurately represented on a coarse
mesh $X_\Delta^d\subset D$, whose maximum mesh spacing 
$\Delta\gg\delta$. Such a
representation can be expressed as, e.g., in terms of polynomial
approximations in the physical space $D$,
\begin{equation} \label{proj_x}
\bar{U}_k(x,t)\simeq
\mathcal{I}_x \bar{U}_k(x,t) = \sum_{l=1}^L \hat{\bar{U}}_{k,l}(t)\phi_k(x),
\end{equation}
where $\{\phi_l(x)\}_{l=1}^L$ are the basis functions in $X_\Delta^d$, and
the expansion coefficients are defined as
$\hat{\bar{U}}_{k,l}=(\bar{U}_k,\phi_l)_\Delta/(\phi_l,\phi_l)_\Delta, 
\forall k$. Here
$(\cdot,\cdot)_\Delta$ is the usual inner product in $X_\Delta^d$ and we have 
assumed, for notational convenience, that the bases are orthogonal.
The completeness of such bases yields
\begin{equation} \label{x_conv}
\Arrowvert\bar{U}_k(x,\cdot)-\mathcal{I}_x \bar{U}_k(x,\cdot)\Arrowvert_
{X_\Delta^d}
=O(L^{-\gamma_{x,k}})\to 0, \qquad L\to\infty,~\gamma_{x,k}>0, \forall k,
\end{equation}
where $\Arrowvert\cdot\Arrowvert_{X_\Delta^d}$ is an appropriate 
norm in $X_\Delta^d$
and $\gamma_{x,k}>0$ quantifies the convergence rate, which depends on the
smoothness of the underlying function $\bar{U}_k$.

Similarly, we can define two operators between $\bar{U}$ and 
$\{\hat{\bar{U}}_l\}$: the {\em restriction operator} $\mathcal{P}_x$
in the physical space $D$ such that
\begin{equation} \label{P_x}
\{\hat{\bar{U}}_{k,l}(\cdot)\}=\mathcal{P}_x
\{\bar{U}_k(x,\cdot)\}, \qquad k=0,\cdots,K \textrm{ and } l=1,\cdots,L,
\end{equation}
and the {\em lifting operator} $\mathcal{Q}_x$ such that
\begin{equation} \label{Q_x}
\mathcal{I}_x
\{\bar{U}_k(x,\cdot)\}=\mathcal{Q}_x\{\hat{\bar{U}}_{k,l}(\cdot)\}, \qquad
k=0,\cdots,K \textrm{ and } l=1,\cdots,L.
\end{equation}
Clearly, we have $\mathcal{Q}_x\mathcal{P}_x=\mathcal{I}_x$ and
$\mathcal{I}_x\mathcal{I}_x=\mathcal{I}_x$.

\subsubsection{Global restriction and lifting} \label{sec:I}

The global {\em restriction operator} $\mathcal{P}$ and 
{\em lifting operator} $\mathcal{Q}$
are thus defined as
\begin{equation} \label{P}
\mathcal{P}=\mathcal{P}_x\mathcal{P}_\omega \textrm{ such that }
\{\hat{\bar{U}}_{k,l}(t)\}=\mathcal{P}u(\omega,x,t);
\end{equation}
and
\begin{equation} \label{Q}
\mathcal{Q}=\mathcal{Q}_\omega\mathcal{Q}_x \textrm{ such that }
\mathcal{I} u(\omega,x,t)=\mathcal{Q}\{\hat{\bar{U}}_{k,l}(t)\},
\end{equation}
where the {\em global approximation operator} $\mathcal{I}$ is defined as
\begin{equation} \label{I}
\mathcal{I}=\mathcal{QP}=
\mathcal{Q}_\omega\mathcal{Q}_x\mathcal{P}_x\mathcal{P}_\omega=
\mathcal{Q}_\omega\mathcal{I}_x\mathcal{P}_\omega.
\end{equation}
Note the above operators are defined for
$\omega=\{\omega_j\}_{j=1}^M\in\Omega$, $k=\{0,\cdots,K\}$, 
$l=\{1,\cdots,L\}$.
A remarkable property of the operator $\mathcal{I}$ is that
\begin{equation} \label{I_conv}
\|u-\mathcal{I}u\|\to 0, \qquad K,L\to\infty,
\end{equation}
where the norm $\|\cdot\|$ is defined in the tensor product space of
$L^2(\Omega)$ and the appropriate space of $D$ that defines the norm
in \eqref{x_conv}. 
Such a property ensures that, by restricting from
the fine-scale solution ensemble to the coarse-scale variables (operator $\mathcal{P}$)
and lifting back to the fine-scale (operator $\mathcal{Q}$), 
we can reconstruct a representative ensemble of fine-scale solutions with
controllable accuracy.

\subsubsection{Implementation of Restriction and Lifting}
\label{sec:details}

In numerical simulations, the ensemble of solutions obtained by the Monte
Carlo method in Section \ref{sec:mcfem}, $\{u(\cdot,\omega_j)\}_{j=1}^M$,
are first restricted in random space by $\mathcal{P_\omega}$
\eqref{P_oper}; the resulting variables $\{\bar{U}_k(x)\}$ on the fine
mesh is further restricted to the coarse mesh by 
$\mathcal{P}_x$ \eqref{P_x}.

The random restriction
$\mathcal{P}_\omega$ is accomplished by \eqref{proj_coef}, where the 
integral is written in a formal way as $u$ and $\xi$ in the integrand
do not in general have the same probability measure. To integrate
\eqref{proj_coef}, we transform both the random variables $u$ and $\xi$
to a uniform variable $\theta\in U(0,1)$ via their cumulative density
function (CDF), i.e.,
$\theta = F(u) = G(\xi)$, where $F$ and $G$ and the CDF of $u$ and
$\xi$, respectively. Hence,
\begin{equation} \label{cdf}
u=F^{-1}(\theta), \qquad \xi=G^{-1}(\theta),
\end{equation}
and \eqref{proj_coef} can be written as
\begin{equation} \label{P_w_cdf}
\bar{U}_k =
\frac{1}{\int_\Omega\Phi_k^2 dP}
\int_\Omega u(\omega,\cdot)\Phi_k(\xi(\omega))dP(\omega)
=\frac{1}{\int_\Omega\Phi_k^2 dP}
\int_0^1 F^{-1}(\theta)\Phi_k(G^{-1}(\theta))d\theta.
\end{equation}
The resulting integral in the bounded domain $(0,1)$ can
be readily integrated using quadrature rules with sufficient
accuracy.
For more details on such integrations and numerical
examples, see \cite{XiuK_SISC02}.

During the lifting procedure \eqref{Q}, the coarse variables on
the coarse mesh are first ``lifted'' to the fine mesh via $\mathcal{Q}_x$
\eqref{Q_x} (effectively, interpolated); the random lifting operation $\mathcal{Q}_\omega$
is then conducted to generate an ensemble of fine-scale solutions
on the fine mesh. 
The spatial lifting $\mathcal{Q}_x$ is accomplished through
\eqref{proj_x}, and the random lifting $\mathcal{Q}_\omega$ by
\eqref{proj_oper}. In \eqref{proj_oper}, an ensemble of realizations
of the random variable $\{\xi(\omega_j)\}_{j=1}^M$ are generated to, 
in turn,  generate $\{u(\omega_j)\}_{j=1}^M$. 
The $\{\xi(\omega_j)\}_{j=1}^M$ is generated via \eqref{cdf} by using the same set of
uniform random variable
$\{\theta(\omega_j)\}_{j=1}^M$ resulted from the CDF of 
$\{u(\omega_j)\}_{j=1}^M$, i.e., $\theta=F(u)$.
By using the same set of $\theta$, the lifted solution 
$\{\mathcal{I}u(\omega_j)\}_{j=1}^M$ will have the same correlation structure as 
$\{u(\omega_j)\}_{j=1}^M$.
This is a key step for efficient numerical computations, as it
prescribes the ``right" correlations between 
a particular realization of the medium and the corresponding
solution for this medium.

\subsection{Coarse-scale integration}

As pointed out in section \ref{sec:I_rand}, 
although the pointwise expansion (\ref{proj_oper}) only utilizes
one-dimensional expansions in random space 
at fixed locations in physical space, it
is unclear, to the authors' best knowledge, 
how the governing equations for the coarse-variables 
$\{\hat{\bar{U}}_{k,l}\}$ should be, i.e. the right-hand-side of
\begin{equation} \label{macro_rhs}
\frac{\partial \hat{\bar{U}}_{k,l}}{\partial t} =
R_{k,l}(\hat{\bar{U}}), \qquad \forall ~k\in [0,K],~ l\in [1,L].
\end{equation}
is unknown.
To circumvent the difficulty, we employ the ``equation-free'' approach
where explicit knowledge of these governing equations is not needed. 
This procedure is as following:
\begin{itemize}

\item Evaluate the coarse variables
  $\{\hat{\bar{U}}_{k,l}(t)\}=\mathcal{P}u(\cdot,t)$ from the
  fine-scale computation, for $t^n\leq t\leq t^n_c=t^n+n_f\delta t$.

\item Approximate the RHS of (\ref{macro_rhs}) at $t=t^n_c$, i.e.
\begin{equation} \label{derivative}
R_{k,l}(t_c^n)=\sum_{j=0}^{n_e}\alpha_j \hat{\bar{U}}_{k,l}(t_j)
=\frac{d\hat{\bar{U}}_{k,l}}{dt}(t^n_c)+O(\delta t^{J_f}),
\qquad k\in[0,K], ~l\in[1,L],
\end{equation}
where $1\leq n_e\leq n_f$, $t_j=t^n_c-j\delta t$,
and $J_f$ denotes the order of the
approximation.
$\{\alpha_j\}_{j=1}^{n_e}$ is a set of consistent 
coefficients such that $\sum\alpha_j v(t_j)=d v/dt(t^n_c) + 
O(\delta t^{J_f})$.

\item Once the RHS of (\ref{macro_rhs}) is estimated numerically, 
\eqref{macro_rhs} is integrated forward in time for one step 
  on a larger time step. For example,
  given coarse time step of size $\Delta t_c = n_c\cdot\delta t$ with
 $n_c\ge 1$, such that $t^{n+1}=t^n+\Delta t_c=t^n+(n_f+n_c)\delta t$,
 the Euler forward integrator takes the form
\begin{equation} \label{coarse_t}
\hat{\bar{U}}_{k,l}^{n+1}=
\hat{\bar{U}}_{k,l}(t^n_c)+\Delta t_c\cdot R_{k,l}(t^n_c)
+O(\Delta t_c^2),
\qquad k\in[0,K], ~l\in[1,L].
\end{equation}
\end{itemize}
The estimation of derivatives \eqref{derivative} may suffer from numerical
oscillations, especially for systems that exhibit noisy behavior at microscopic
scales, e.g., molecular dynamics. In this case, certain smoothing techniques such
as a least-sqaure fit may be used to alleviate the problem \cite{GearLK_PLA03}.
(Such is not the case in this paper.) The integration of coarse variables
\eqref{coarse_t} is a simple Euler forward scheme, which can lead to
relatively large errors when $\Delta t_c$ is large.
Other integration
schemes, such as higher-order single step methods or multi-step
methods, can be used for their improved accuracy and/or
better stability properties \cite{GearK_SISC03, GearK_JCP03}.
A complete analysis (stability, accuracy, etc.) 
of the equation-free method to stochastic equations
is still lacking, and is beyond the scope of the current paper. For an error
analysis of equation-free method to a deterministic system (flow simulation),
see \cite{SirisupXKK_JCP05}.

\subsection{Computational Complexity}

Let us denote by $N_f\sim O(\delta^{-d})$ the number of
degrees-of-freedom (DOF) of
the fine mesh $X_\delta^d$, and $N_c\sim O(\Delta^{-d})$ the DOF of the coarse
mesh $X_\Delta^d$. During the coarse integration step, the 
fine-scale computation by $M$ realizations of
Monte Carlo simulations is effectively reduced to a problem of
the evolution of $(K+1)$ local polynomial expansion coefficients
(\ref{proj_oper}) on
the coarse mesh obtained by the global restriction operator $\mathcal{P}=
\mathcal{P}_x\mathcal{P}_\omega$. The reduction in computational complexity
can be illustrated as
\begin{equation} \label{map}
\mathbb{R}^{M\times N_f} \stackrel{\mathcal{P}_\omega}{\longrightarrow}
\mathbb{R}^{(K+1)\times N_f} \stackrel{\mathcal{P}_x}{\longrightarrow}
\mathbb{R}^{(K+1)\times N_c}.
\end{equation}
Thus, to march problem (\ref{multi_eqn}) over a global time step
$\Delta t=(n_f+n_c)\delta t$, the particular equation-free 
projective integration algorithm, which consists
of $n_f$ steps of fine-scale computations and one-step of coarse 
integration, 
needs to solve a problem of complexity
\begin{equation} \label{complex}
C_c \sim n_f\times\mathbb{R}^{M\times N_f}+\mathbb{R}^{(K+1)\times
  N_c}.
\end{equation}
On the other hand, the complexity of the full-scale MCFEM over the
same time interval $\Delta t$ is
\begin{equation} \label{complex_f}
C_f \sim (n_f+n_c)\times\mathbb{R}^{M\times N_f}
\end{equation}
For the multiscale problem considered in this paper, we have
$N_f\gg N_c$, $O(1)\sim K\ll M$, and $n_c\gg n_f$. Thus, roughly
speaking, $C_f/C_c\sim (1+n_c/n_f)\gg1$. 
We remark that such an
estimate is rather crude and the actual
computational efficency is problem dependent.

\section{Numerical Results} \label{sec:results}

In this section we present numerical results on (\ref{multi_eqn})
in one spatial dimension $x\in[0,1]$. 
The restriction and lifting operations
in physical spaces, as shown
in Section \ref{sec:I_phys},
come from standard approximation theory, and here we focus on
the properties of the {\em random restriction} and the {\em random
lifting}. 
We remark that the method extends 
trivially to multi-dimensional physical spaces.
We assume that $\kappa(\omega,x)$ is a Gaussian random field with
unit mean value, i.e., $\mathbb{E}u(\omega,x)=1$, and employ the
Hermite polynomials in random space
to represent the random field in \eqref{proj_oper}.
We employ the Gaussian random field model because of its ease of
numerical generation. (Generation of non-Gaussian random fields
is still an active research area.)
From a mathematical point of view, 
Gaussian models are inappropriate for the diffusivity fields as they allow
negative values with non-zero probability, and thus violate the
uniform coercivity assumption (\ref{coercivity}). In practice, however,
such negative values are rare, especially when the variance is small,
and we can neglect such negative values if they occur in the random 
realizations.
(In the numerical examples below, negative values never
appear.) 

We assume that the Gaussian random field has an 
exponential covariance function, i.e.
$C_{\kappa}(x,y)=\exp(|x-y|)/l_\kappa$. Such a correlation function
can be generated from a first-order Markov process, and has been used
extensively in the literature. 
All MCFEM computations are conducted by a linear
finite element method with Euler forward integration, and $M=1,000$
realizations are used. The forcing term in (\ref{multi_eqn}) is fixed
at $f(x)\equiv -2$, and zero Dirichlet boundary conditions are imposed
at $x=0$ and $x=1$. The restriction in physical space (\ref{proj_x}) is
conducted on a set of Jacobi polynomial basis (see \cite{KarniadakisS99},
Chapter 2).

\subsection{Accuracy}

The complete error analysis of the present multiscale method
remains an open issue.  In this section, we conduct numerical experiments
to document the error convergence.
In the first example, we set $l_\kappa=0.1$ and use 40 linear elements
($\delta=0.025$).
Fifth-order Hermite expansion ($K=5$) is used for the random
restriction (\ref{proj_oper}), and fifth-order Jacobi
basis ($L=5$) for the spatial restriction (\ref{proj_x}).
The time step for the
fine-scale computation is $\delta t=0.001$ and it is conducted for $n_f=20$
steps. 
The coarse integration has a time step $\Delta t_c=0.08$ (i.e. $n_c=80$),
so that the global
time step is $\Delta t=(n_f+n_c)\Delta t_f=0.1$. 
For this illustrative example, the computational speed-up is modest
($4\sim 5$ times), because the separation of scales is modest.

\begin{figure}[htbp]
   \centerline{
   \psfig{file=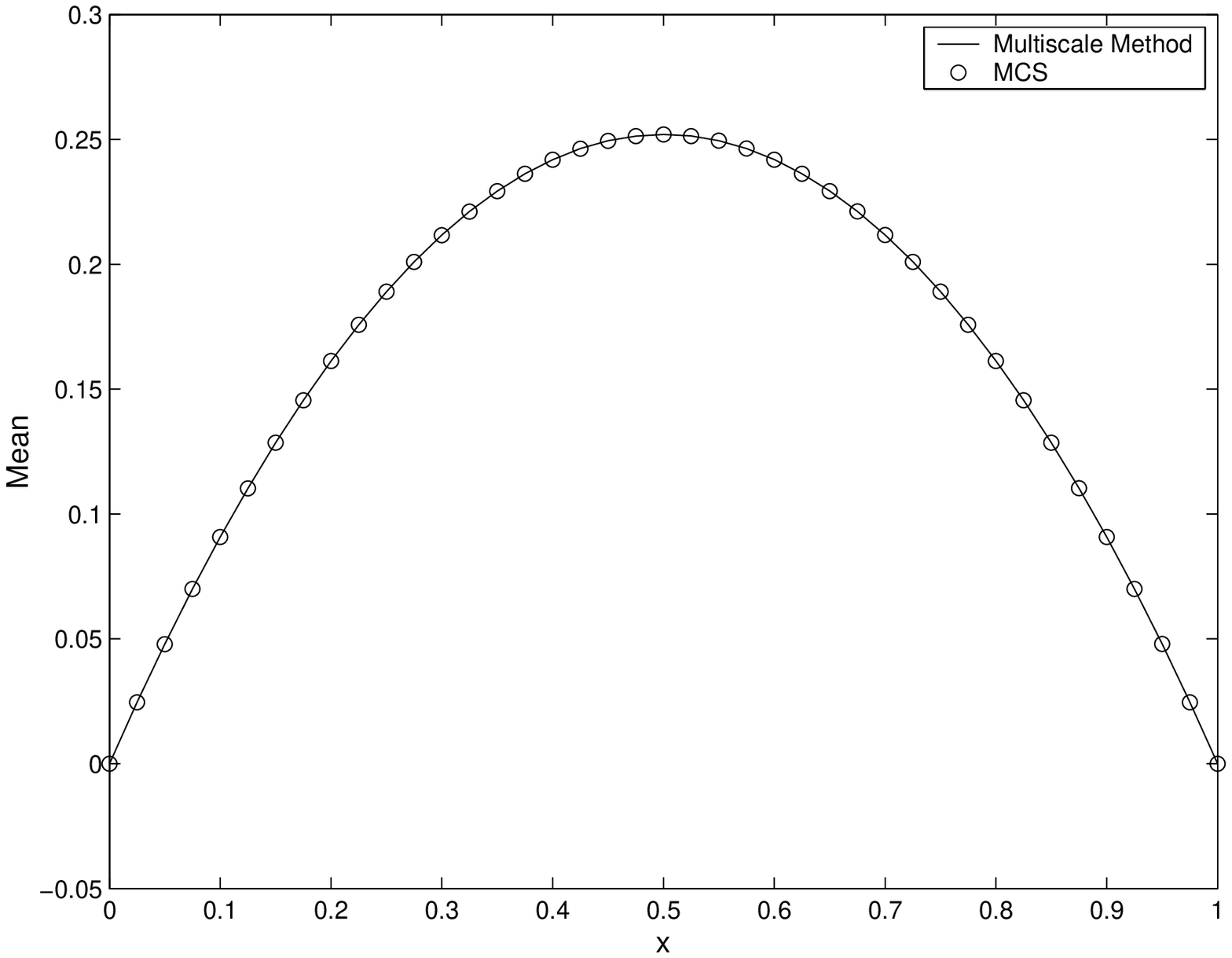,width=6cm}\qquad
   \psfig{file=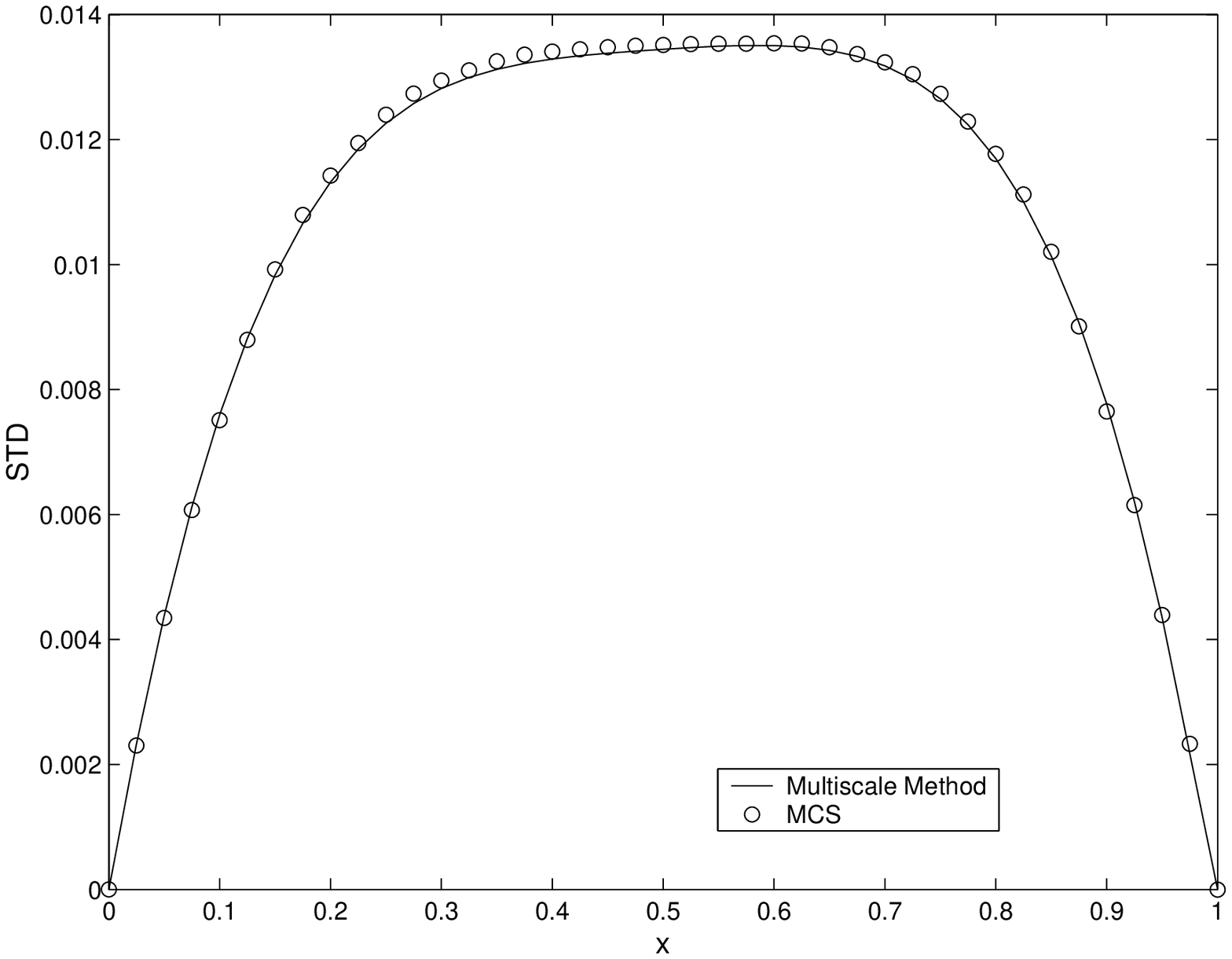,width=6cm}}
\caption{Solution profile at $T=1$.
         Left: mean solution,
         Right: standard deviation.
        }
\label{fig:run1}
\end{figure}

Figure \ref{fig:run1} shows the stochastic solution profile (mean and standard
deviation) at time $T=1$. Good agreements are obtained between
the full-scale MC simulation and the equation-free multiscale method.

\begin{figure}[htbp]
   \centerline{
   \psfig{file=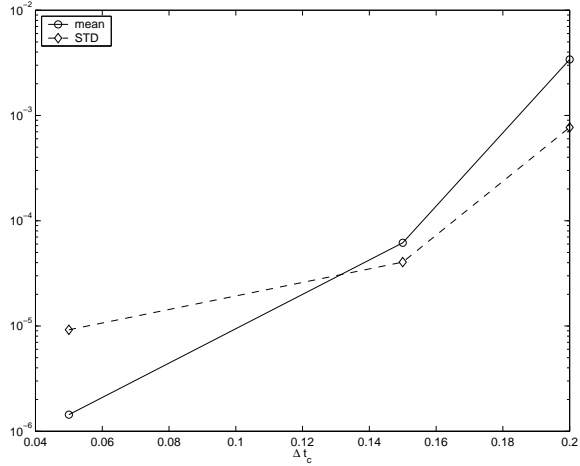,width=3in}}
\caption{Errors in mean and standard deviation (STD) versus the step size
         of the coarse integration.}
\label{fig:errDt}
\end{figure}

To study the error contributions from different factors, we
conducted a series of computations with varying parameters.
In Figure \ref{fig:errDt} the $L^\infty$ 
errors in mean and standard deviation (STD) are shown. These
computations have fixed values of
$\Delta t_f=0.05, K=3, L=4$ and varying time steps $\Delta t_c$ 
of the coarse integration. 
We observe that the errors decrease as the size of time steps for the
coarse integration decreases.

\begin{figure}[htbp]
   \centerline{
   \psfig{file=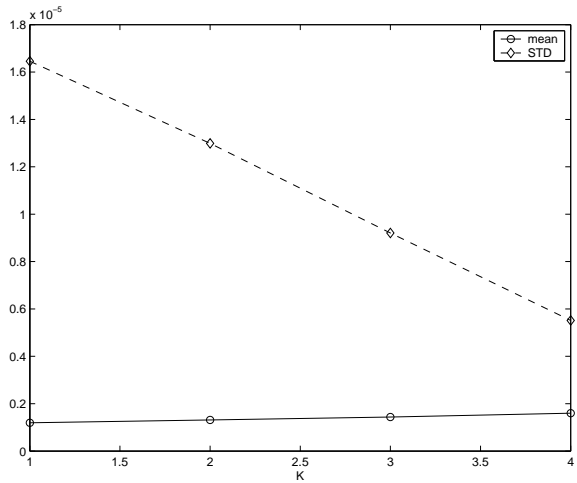,width=3in}}
\caption{Errors in mean and standard deviation (STD) versus the order
         of Hermite approximation of the coarse variables.}
\label{fig:errK}
\end{figure}

To examine the error convergence with respect to the orders of 
approximation in random space (parameter $K$) and physical space 
(parameter $L$), we fix
the time steps of integrations $\Delta t_f=\Delta t_c=0.05$. The size
of coarse integration $\Delta t_c$ is sufficiently small such that the
the temporal errors are subdominant (e.g., $O(10^{-6})$ for the mean as shown
from Fig. \ref{fig:errDt}). In Fig. \ref{fig:errK}, the errors with
increasing order of the Hermite approximations $(K)$ 
in the random restriction $\mathcal{P}_\omega$ are shown.
Fourth-order ($L=4$) Jacobi basis is used in the physical
space, so that the errors from spatial restriction $\mathcal{P}_x$ are
subdominant. It can be seen that as the order $K$ of the random restriction
increases, the errors in standard deviation decrease as expected.
The mean solution is well-resolved by even the first-order Hermite expansion
($K=1$) and its errors remains at the $O(10^{-6})$ level.

\begin{figure}[htbp]
   \centerline{
   \psfig{file=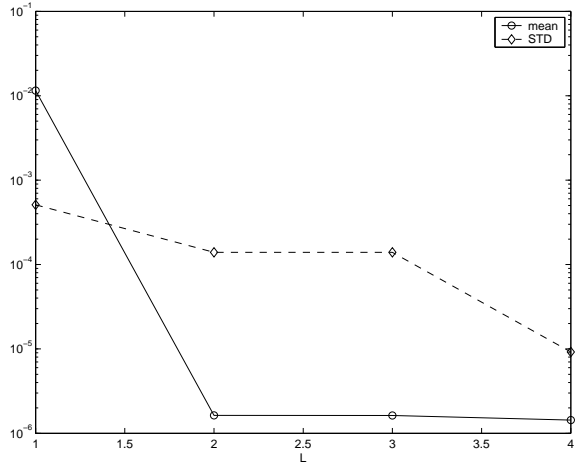,width=3in}}
\caption{Errors in mean and standard deviation (STD) versus the order
         of spatial representation of the coarse variables.}
\label{fig:errL}
\end{figure}

We then fix the order of the Hermite approximation in the 
random space at third-order,
i.e., $K=3$, and vary $L$. Again, $\Delta t_c=\Delta t_f=0.05$ is sufficiently
small. In Fig. \ref{fig:errL}, it can be seen that the error in the mean
solution
quickly reaches a saturation level of $O(10^{-6})$ at second-order $L=2$.
This is consistent with the result in Fig. \ref{fig:errK}. The error
in STD keeps decreasing with increasing order of $L$.

\begin{figure}[htbp]
   \centerline{
   \psfig{file=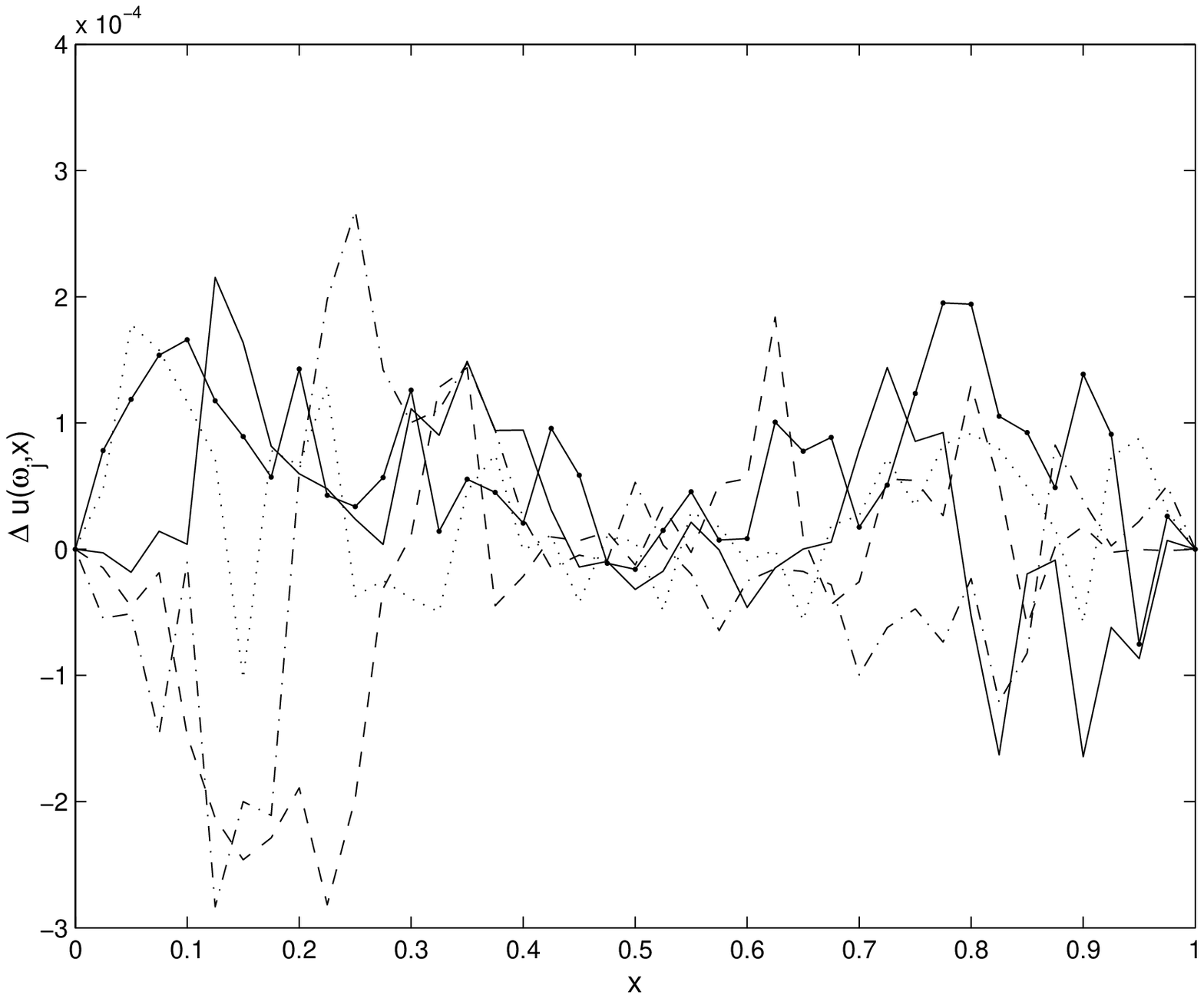,width=6cm}\qquad
   \psfig{file=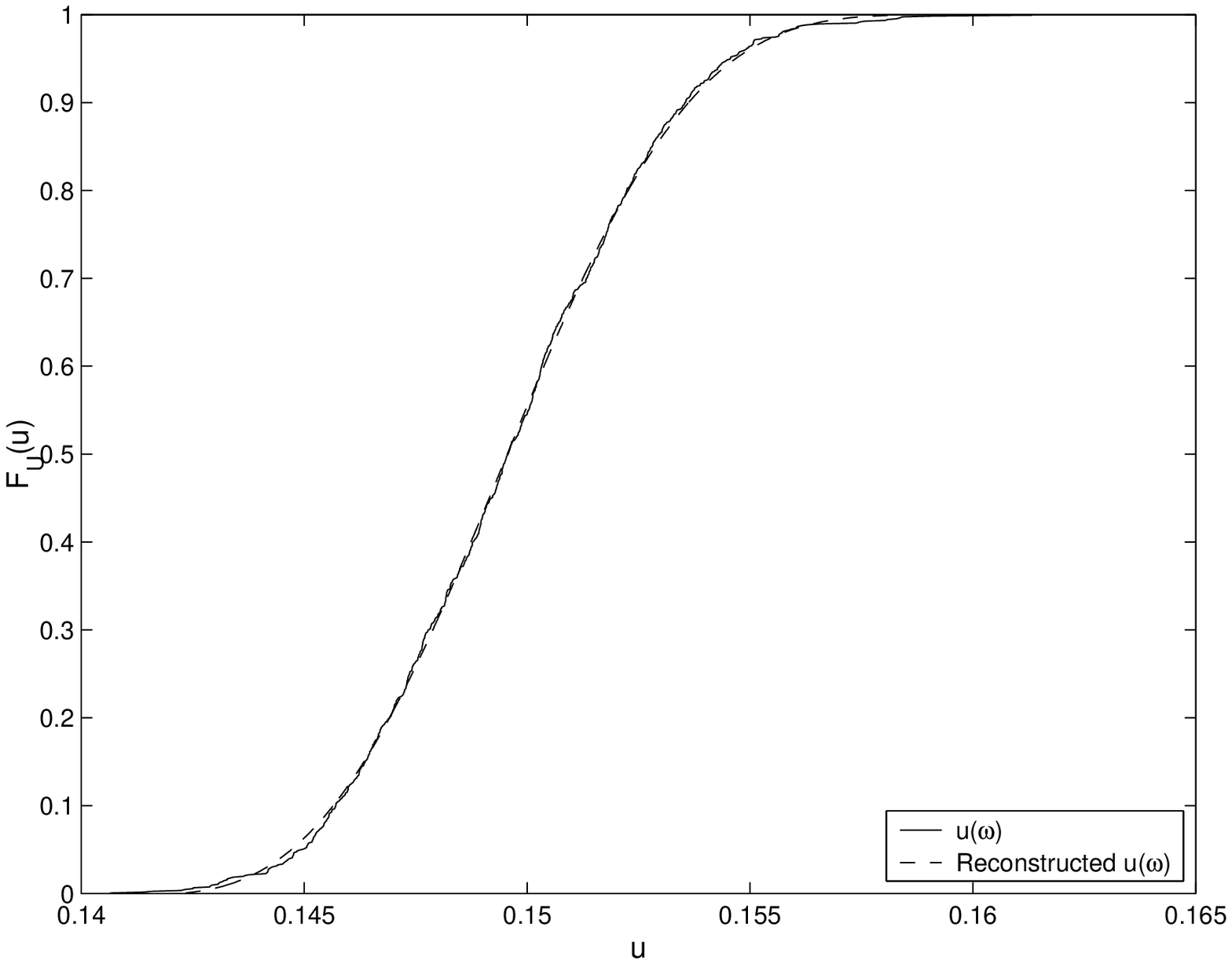,width=6cm}}
\caption{Comparison of random solution ($u$) and its lifting
         ($\mathcal{I}u$) at $t=0.1$.
         Left: profiles of $u(\omega_j,x)-\mathcal{I}u(\omega_j,x)$ for some
               realizations of $\omega_j$,
         Right: CDFs of $u$ and $\mathcal{I}u$ at  $x=0.5$.
        }
\label{fig:proj1}
\end{figure}

As discussed in Section \ref{sec:I}, the current implementations of the
restriction operator $\mathcal{P}$ and the lifting operator 
$\mathcal{Q}$ allow us to
reconstruct representative ensembles of fine-scale solutions 
by restricting them to the coarse
scale first and then lifting back to the fine scale, i.e. $\mathcal
{I}=\mathcal{QP}$ is an approximation operator.
To examine the properties of the global approximation operator
$\mathcal{I}$ (\ref{I}),
we plot $\Delta u(\omega,x,t)=u(\omega,x,t)-\mathcal{I}u(\omega,x,t)$ 
at an arbitrary chosen time $t=0.1$. On the left of Figure \ref{fig:proj1},
several realizations of such errors
(randomly chosen from the $M=1,000$ realizations) are shown. It can be
seen that the errors are bounded within a small range of the
same order of the spatial error ($O(10^{-3})$). 
The CDF of the random solution $u$ at
the center of the domain ($x=0.5$) is shown on the right of Fig.
\ref{fig:proj1}. Again we see excellent agreement between the probability
distribution of the target
$u$ and that of its lifting $\mathcal{I}u$. 

\begin{figure}[htbp]
   \centerline{
   \psfig{file=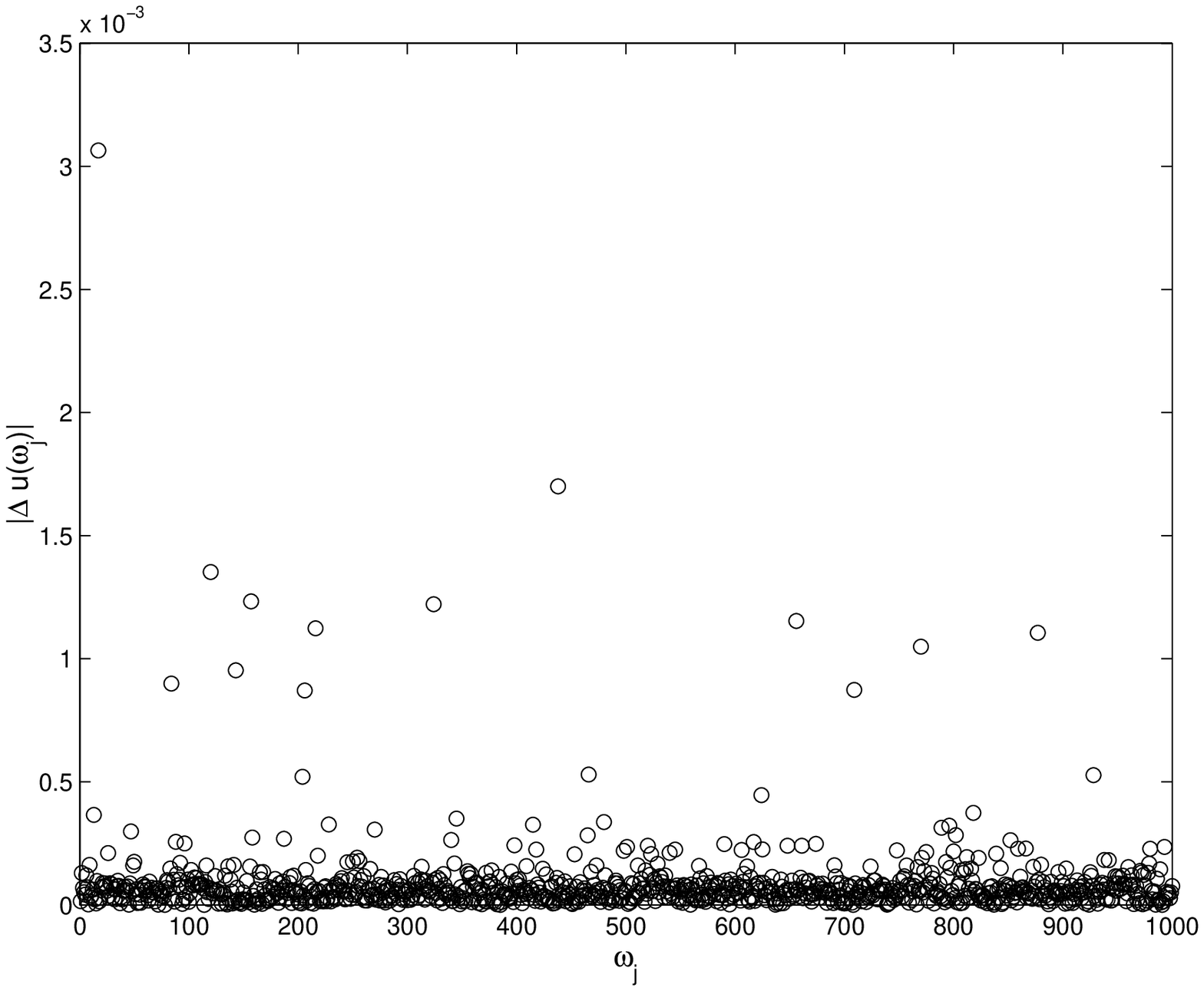,width=6cm}\qquad
   \psfig{file=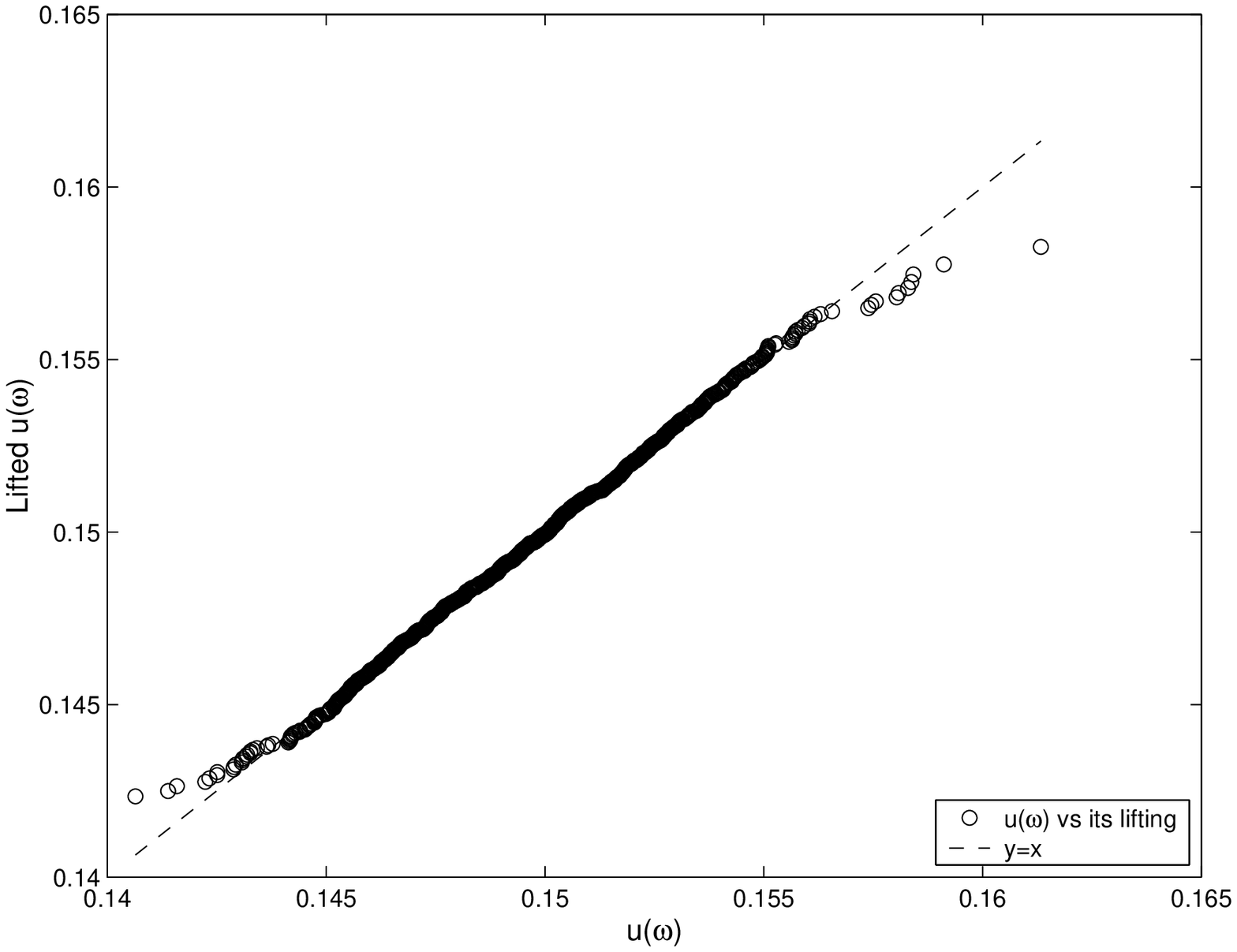,width=6cm}}
\caption{Comparison of random solution ($u$) and its lifting
         ($\mathcal{I}u$) at $t=0.1, x=0.5$.
         Left: point-wise error of $|u(\omega_j)-\mathcal{I}u(\omega_j)|$ for
                all $\omega_j, j=1,\dots,1,000$.
         Right: $u(\omega_j)$ vs. $\mathcal{I}u(\omega_j)$.
        }
\label{fig:proj2}
\end{figure}

\begin{figure}[htbp]
   \centerline{
   \psfig{file=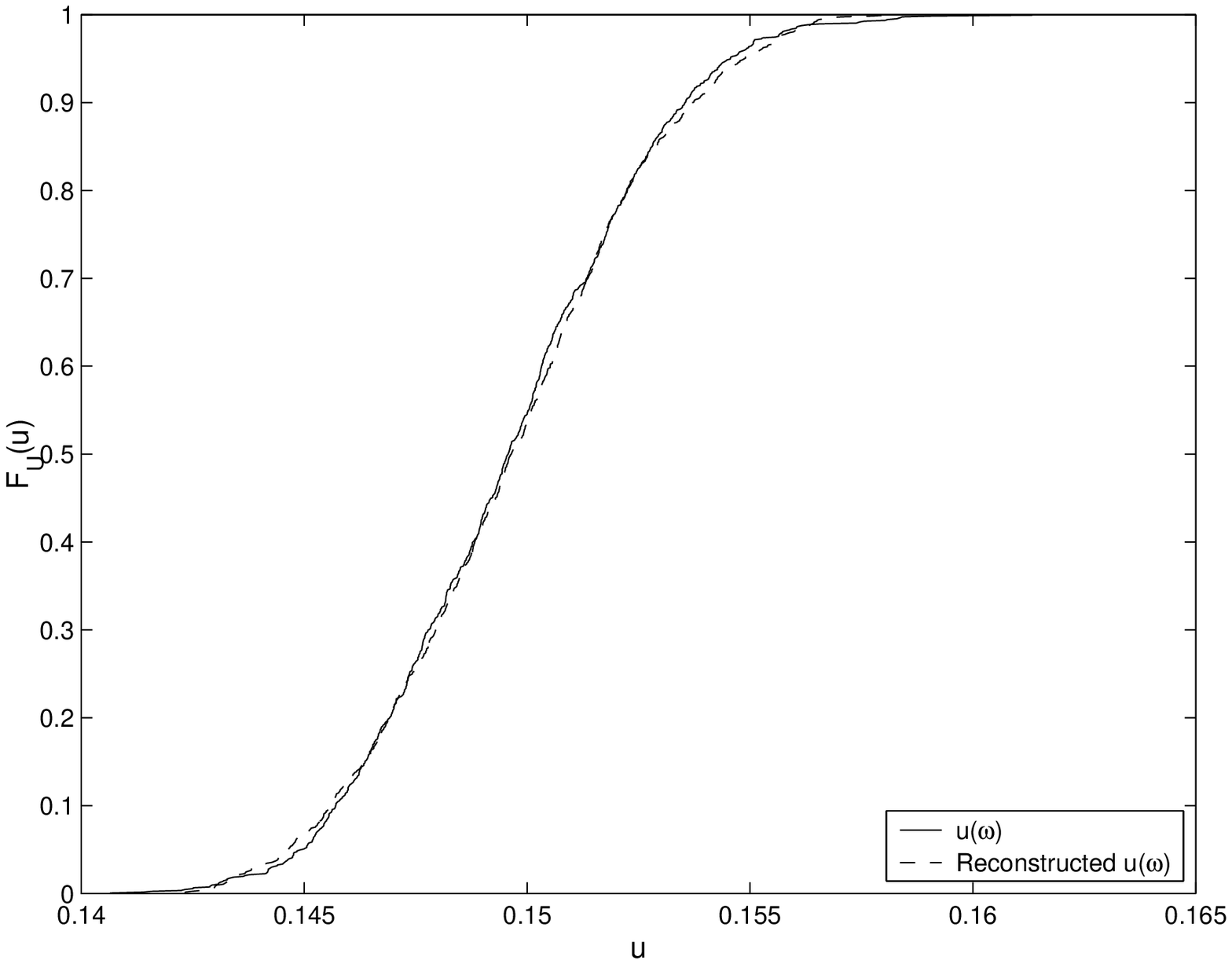,width=6cm}\qquad
   \psfig{file=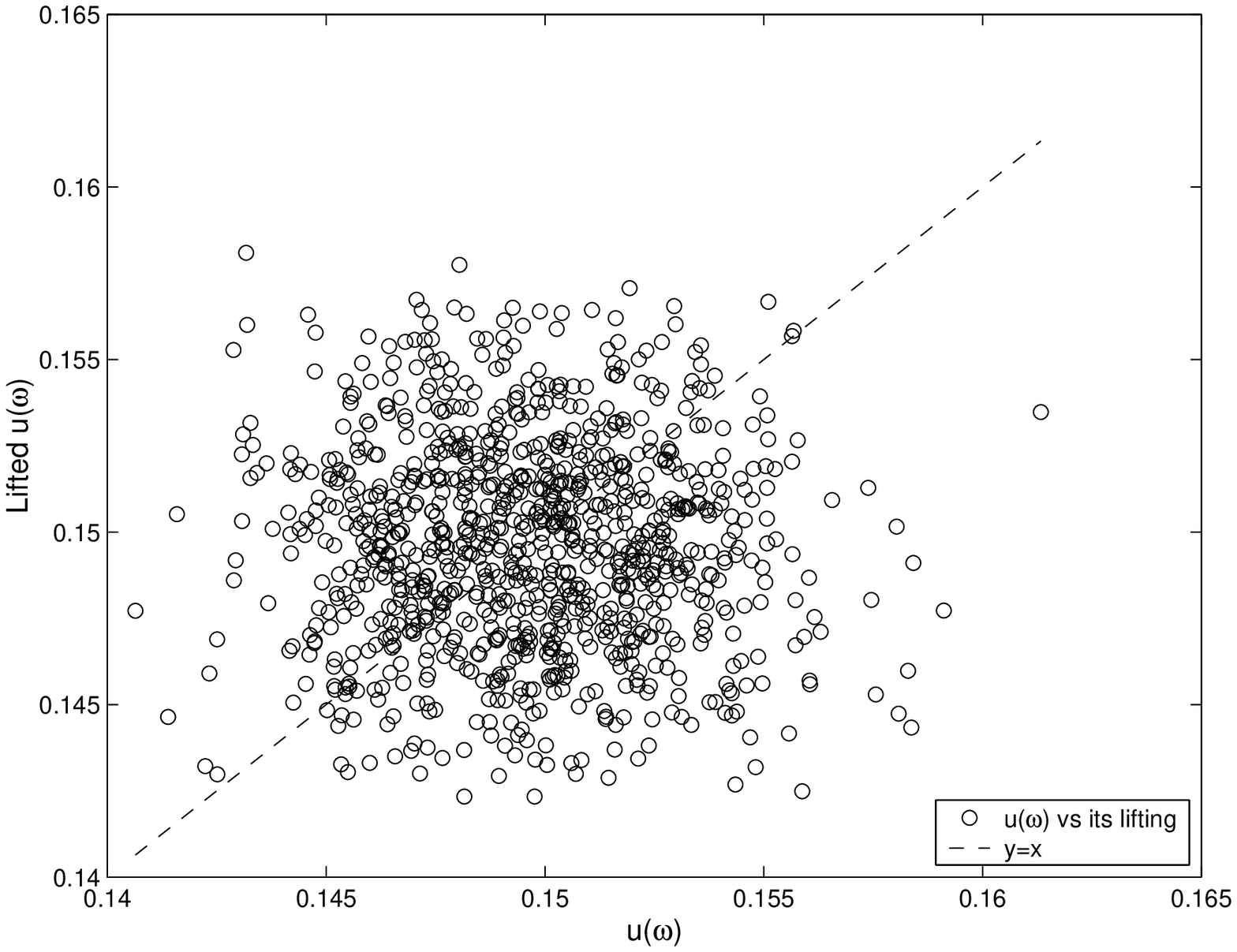,width=6cm}}
\caption{Comparison of random solution ($u$) and its lifting
         ($\mathcal{I}u$) without correct correlation
          structure at $t=0.1, x=0.5$.
         Left: CDF of $u$ and $\mathcal{I}u$;
         Right: $u(\omega_j)$ vs. $\mathcal{I}u(\omega_j)$ for $j=1,\dots,
                 1,000$.
        }
\label{fig:proj_bad}
\end{figure}

To further examine the error of $\mathcal{I}$, we plot on the left of
Fig. \ref{fig:proj2} the point-wise error of $\Delta u(\omega,x,t)$ at
$x=0.5, t=0.1$ for all realizations $\omega_j, j\in[1,1,000]$. 
We observe that,
except at a few discrete points, which belong to a set with {\em
arguably} zero measure,
the errors are bounded in a very small interval of order $O(10^{-3})$. 
On the right of 
Fig. \ref{fig:proj2}, we plot the path-wise correspondence of
$u(\omega_j)$ vs. $\mathcal{I}u(\omega_j)$ for all $j=1,\cdots,1,000$.
It is seen that the data collapse on $y=x$ where the exact correspondence
should be. 
These results confirm that our $\mathcal{I}$ is indeed an approximation
operator which allows us to 
reconstruct the fine-scale solution ensemble {\em with built in desired correlation
structures} from the computed coarse solutions. 
To achieve this, it is
important that one uses, during the lifting step \eqref{proj_oper},
the same uniform random variable values 
determined by the numerical solutions at the restriction stage as described in 
Section \ref{sec:details}.  
If, however, a random variable $\xi$ is
used without maintaining the correct correlation structure
between the medium realization and the solution in this medium, the lifted
solutions will not be properly correlated to the true solutions, even if the
ensemble is constructed to maintain the same distribution. 
Figure \ref{fig:proj_bad}
shows such an example. 
Again, this is the numerical solution at $t=0.1$
and we plot the realizations at $x=0.5$. Here the numerical solutions
are lifted by using arbitrarily generated 
Gaussian random variables $\xi$ in \eqref{proj_oper}. We observe that
although the solution has the same distribution as the MCS solutions
(Fig. \ref{fig:proj_bad}, left), it is completely uncorrelated
to the true solution as shown on the right of Fig. \ref{fig:proj_bad}.
The ability to maintain good correlation structure in the
lifting procedure is a distinctive feature of our method.
This is different from the conventional lifting
procedures, whose reconstructed solutions are rather arbitrary (to a certain
degree) and
a certain constraining procedure or a relaxation integration is needed
to ``heal'' the lifted solution ensemble \cite{GearK_JSC04,GearKKZ_SIADS04}.

\subsection{Efficiency}

\begin{figure}[htbp]
   \centerline{
   \psfig{file=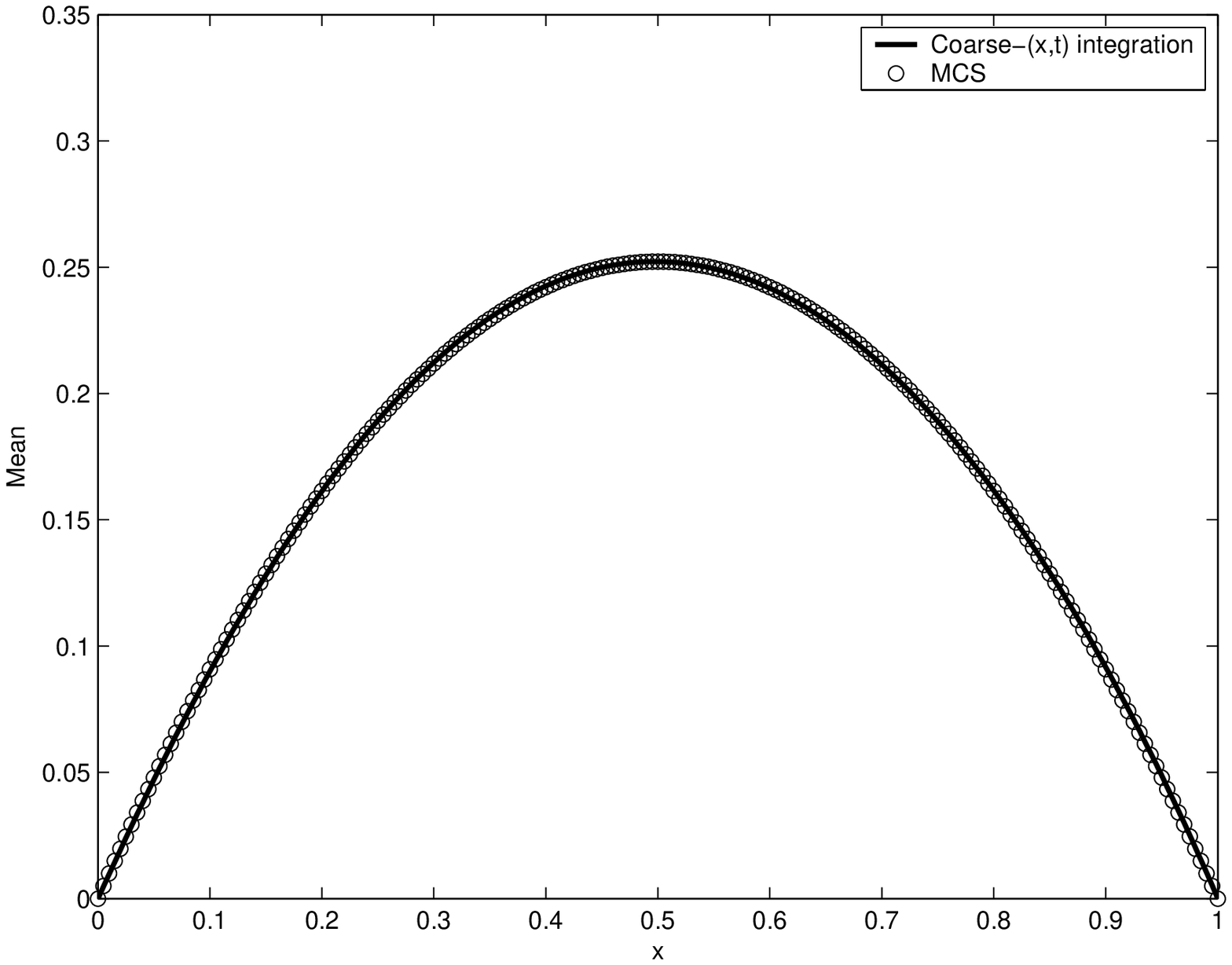,width=6cm}\qquad
   \psfig{file=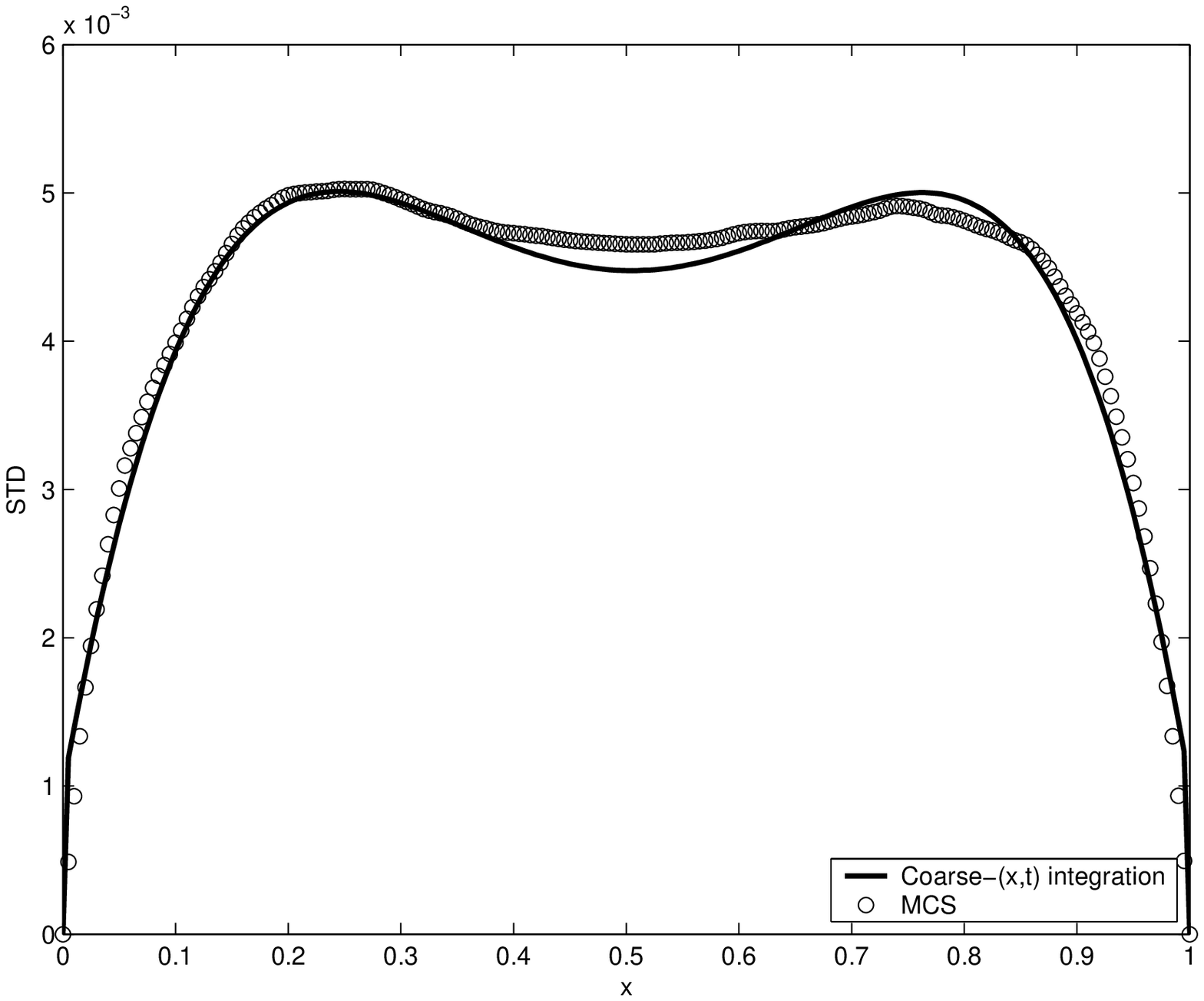,width=6cm}}
\caption{Solution profile at $T=1$.
         Left: mean solution,
         Right: standard deviation.
        }
\label{fig:run2}
\end{figure}

In the second example, we prescribe a diffusivity field with relatively 
small correlation
length $l_\kappa=0.01$. 
We employ $1,000$ linear elements to resolve
the small spatial scales ($\delta=0.001$), 
and this results in a time scale of
$O(10^{-6})$. Thus, we set the time step of the fine-scale computation
at $\delta t=10^{-6}$.
The number of
fine-scale computations within each global time step is $n_f=1000$,
and the time step of the coarse integration is chosen as
$\Delta t_c=n_c\delta t$ with $n_c=49,000$. Thus, 
the global time step is $\Delta t=(n_f+n_c)\delta t=0.05$. 
The polynomial orders are set at 
$K=5$ and $L=5$ for the restrictions in the random space and the physical
space, respectively.
For this application,
we achieve computational speed-up of $\sim 50$,
compared to the full-scale MC simulation.
The random solution profiles in Figure \ref{fig:run2} show again good
agreement between the full-scale MCFEM and the multiscale
computation. 
We remark that the quantification of the 
computational speed-up
is problem dependent. In the diffusion problems considered here, such
speed-up is smaller at the beginning of the computation due to the fast
evolution of the solution. However, the speed-up is significantly larger
once the initial transient is passed.

\section{Summary} \label{sec:summary}

In this paper we present an equation-free multiscale algorithm for 
integrating unsteady
diffusion problems in a random medium with small-scale spatial structures,
e.g., short correlation length.
The method is based on the assumption that although the individual
realizations of the random solutions are characterized by small scales, 
their ensemble averages are much smoother, characterized
by larger scales. 
This motivates the use a set of coarse-grained observation variables,
which are based on the pointwise polynomial
approximations of the fine-scale random solutions. 
Such coarse-grained variables are approximated accurately on a coarse mesh,
and integrated in time with large time steps. 
An equation-free approach
is employed for the coarse integration, as the explicit knowledge of
the governing equations of the coarse variables is unavailable in closed form.
Details of the multiscale method is presented, and
its accuracy and efficiency are documented by numerical examples. 
In particular,
we demonstrate that the present constructions of the restriction and lifting
operators allow us to successfully reconstruct representative fine-scale solutions based
only on the knowledge of the coarse solutions. 
Future work will
include a complete analysis of error estimates of the current
method and applications to more complicated systems.

It is interesting that, in this approach, one creates a ``wrapper" around
an existing direct detailed solver, using it as a black box; traditional
stochastic Galerkin algorithms would require the writing of new code to solve
the coupled system of equations in both physical and random space.
Our approach can thus be considered as a ``nonintrusive" one -- the solution
in both physical and random space is solved using an existing
legacy code through a wrapper, and sidestepping the effort of new 
code development and validation.
In this paper the only
numerical task we demonstrated in the equation-free context was
{\it temporal integration}.
Other tasks enabled through
matrix-free iterative linear algebra (e.g. Newton-Krylov
GMRES based fixed point solvers, Arnoldi-type eigensolvers) 
naturally fit in the equation-free framework; and many more
can be performed on the explicitly unavailable coarse-grained
equation: steady state and bifurcation computations, stability computations,
equation-free optimization and even dynamic renormalization.
Finally, the smoothness {\it in space} of the coarse-grained variables
can be exploited (via the so-called gap-tooth and patch dynamics equation
free schemes) to perform the fully resolved fine-scale computations
not only for short times, but also for only parts of the physical domain 
of interest \cite{KevrekidisGHKRT_CMS03,SamaeyRK_SIMMS04,GearLK_PLA03}.
%
%
Finally, although in this paper we have designed direct numerical experiments to {\it solve} 
the coarse-grained equation which is hypothesized to exist and close,
it is worth noting that it is also possible to design
direct numerical experiments to {\it test this hypothesis} (see \cite{LiKGK_SIMMS03}).
These tasks, and the conditions under which they can be successfully performed
in an equation-free framework and 
accelerate random computations, is the subject of ongoing research.

\section*{Acknowledgments}
We would like to thank Dr. Ivo Babu$\breve{\textrm{s}}$ka and Dr. Raul Tempone
of ICES of the University of Texas, Austin, for useful discussions. 
This work was partially supported through DARPA, AFOSR and an NSF/ITR grant.

\bibliographystyle{siam}
\bibliography{multiscale,random}

\begin{thebibliography}{10}

\bibitem{BensoussanLP78}
{\sc G.~P. A.~Bensoussan, J.-L.~Lions}, {\em Asymptotic analysis for periodic
  structures}, North-Holland, Amsterdam, 1978.

\bibitem{BabuskaC_CMAME02}
{\sc I.~Babu$\breve{\textrm{s}}$ka and P.~Chatzipantelidis}, {\em On solving
  elliptic stochastic partial differential equations}, Comput. Methods Appl.
  Mech. Engrg., 191 (2002), pp.~4093--4122.

\bibitem{BabuskaTZ02}
{\sc I.~Babu$\breve{\textrm{s}}$ka, R.~Tempone, and G.~Zouraris}, {\em Galerkin
  finite element approximations of stochastic elliptic differential equations},
  SIAM J. Numer. Anal., 42 (2004), pp.~800--825.

\bibitem{CioranescuD99}
{\sc D.~Cioranescu and P.~Donato}, {\em An introduction to homogenization},
  Oxford University Press, Oxford, 1999.

\bibitem{Dagan89}
{\sc G.~Dagan}, {\em Flow and transport in porous formations}, Springer-Verlag,
  Heidelberg, Berlin, New York, 1989.

\bibitem{DebBO01}
{\sc M.~Deb, I.~Babu$\breve{\textrm{s}}$ka, and J.~Oden}, {\em Solution of
  stochastic partial differential equations using {Galerkin} finite element
  techniques}, Comput. Methods Appl. Mech. Engrg., 190 (2001), pp.~6359--6372.

\bibitem{E_CPA92}
{\sc W.~E}, {\em Homogenization of linear and nonlinear transport equations},
  Comm. Pure Appl. Math., XLV (1992), pp.~301--326.

\bibitem{Fishman96}
{\sc G.~Fishman}, {\em {Monte Carlo: Concepts, Algorithms, and Applications}},
  Springer-Verlag New York, Inc., 1996.

\bibitem{Freeze_WRR75}
{\sc R.~Freeze}, {\em A stochastic-concepcual analysis of one-dimensional
  groundwater flow in nonuniform homogeneous media}, Water Resour. Res., 11
  (1975), pp.~725--741.

\bibitem{GearKKZ_SIADS04}
{\sc C.~Gear, T.~Kaper, I.~Kevrekidis, and A.~Zagaris}, {\em Projecting on a
  slow manifold: singularly perturbed systems and legacy codes}, SIAM J. Appl.
  Dyn. Sys., submitted (2004).
\newblock (original version can be found as Physics/0405074 at arXiv.org).

\bibitem{GearK_SISC03}
{\sc C.~Gear and I.~Kevrekidis}, {\em Projective methods for stiff differential
  equations: problems with gaps in their eigenvalue spectrum}, SIAM J. Sci.
  Comput., 24 (2003), pp.~1091--1106.
\newblock (also NEC Technical Report NECI-TR 2001-029, can be obtained as
  http://www.neci.nj.nec.com/homepages/cwg/projective.pdf).

\bibitem{GearK_JCP03}
\leavevmode\vrule height 2pt depth -1.6pt width 23pt, {\em Telescopic
  projective methods for parabolic differential equations}, J. Comput. Phys.,
  187 (2003), pp.~95--109.
\newblock (also a NEC Technical Report, November 2001, can be obtained as
  http://www.neci.nj.nec.com/homepages/cwg/itproje.pdf).

\bibitem{GearK_JSC04}
\leavevmode\vrule height 2pt depth -1.6pt width 23pt, {\em Constraint-defined
  manifolds: a legacy-code approach to low-dimensional computations}, J. Sci.
  Comput., submitted (2004).

\bibitem{GearLK_PLA03}
{\sc C.~Gear, J.~Li, and I.~Kevrekidis}, {\em The gap-tooth method in particle
  simulations}, Phys. Lett. A, 316 (2003), pp.~190--195.

\bibitem{GhanemS91}
{\sc R.~Ghanem and P.~Spanos}, {\em Stochastic Finite Elements: a Spectral
  Approach}, Springer-Verlag, 1991.

\bibitem{HouW_JCP97}
{\sc T.~Hou and X.-H. Wu}, {\em A multiscale finite element method for elliptic
  problems in composite materials and porous media}, J. Comput. Phys., 134
  (1997), pp.~169--189.

\bibitem{HouWC_MC99}
{\sc T.~Hou, X.-H. Wu, and Z.~Cai}, {\em Convergence of a multi-scale finite
  lement method for elliptic problems with rapidly oscillating coefficients},
  Math. Comput., 68 (1999), pp.~913--943.

\bibitem{HughesFMQ_CMAME98}
{\sc T.~Hughes, G.~R. Feijoo, L.~Mazzei, and J.-B. Quincy}, {\em The
  variational multiscale method -- a paradigm for computational mechanics},
  Comput. Meth. Appl. Math. Engrg., 166 (1998), pp.~3--24.

\bibitem{KarniadakisS99}
{\sc G.~Karniadakis and S.~Sherwin}, {\em Spectral/hp Element Methods for CFD},
  Oxford University Press, 1999.

\bibitem{ShortManifesto}
{\sc I.~Kevrekidis, C.~Gear, and G.~Hummer}, {\em Equation-free: the
  computer-assisted analysis of complex, multiscale systems}, A.I.Ch.E Journal,
  50 (2004), pp.~1346--354.

\bibitem{KevrekidisGHKRT_CMS03}
{\sc I.~Kevrekidis, C.~Gear, J.~Hyman, P.~Kevrekidis, O.~Runborg, and
  C.~Theodoropoulos}, {\em Equation-free coarse-grained multiscale computation:
  enabling microscopic simulators to perform system-level analysis}, Comm.
  Math. Sci., 1 (2003), pp.~715--762.
\newblock (original version can be found as physics/0209043 at arXiv.org).

\bibitem{KleiberH92}
{\sc M.~Kleiber and T.~Hien}, {\em The stochastic finite element method}, John
  Wiley \& Sons Ltd, 1992.

\bibitem{MaitreKNG_JCP04}
{\sc O.~{Le Maitre}, O.~Knio, H.~Najm, and R.~Ghanem}, {\em Uncertainty
  propagation using {Wiener-Haar} expansions}, J. Comput. Phys., 197 (2004),
  pp.~28--57.

\bibitem{MaitreNGK_JCP04}
{\sc O.~{Le Maitre}, H.~Najm, R.~Ghanem, and O.~Knio}, {\em Multi-resolution
  analysis of {Wiener-type} uncertainty propagation schemes}, J. Comput. Phys.,
  197 (2004), pp.~502--531.

\bibitem{LiKGK_SIMMS03}
{\sc J.~Li, P.~Kevrekidis, C.~Gear, and I.~Kevrekidis}, {\em Deciding the
  nature of the coarse integration through microscopic simulations: the
  baby-bathwater scheme}, SIAM J. Multiscale Model. Simulation, 1 (2003),
  pp.~391--407.

\bibitem{LiuBM86}
{\sc W.~Liu, G.~Besterfield, and A.~Mani}, {\em {Probabilistic finite elements
  in nonlinear structural dynamics}}, Comp. Meth. Appl. Mech. Eng., 56 (1986),
  pp.~61--81.

\bibitem{Loeve77}
{\sc M.~Lo\`{e}ve}, {\em Probability Theory, Fourth edition}, Springer-Verlag,
  1977.

\bibitem{MakeevMK_JChP02}
{\sc A.~Makeev, D.~Maroudas, and I.~Kevrekidis}, {\em Coarse stability and
  bifurcation analysis using stochastic simulators: kinetic {Monte Carlo}
  examples}, J. Chem. Phys., 116 (2002), pp.~10083--10091.

\bibitem{MakeevMPK_JChP02}
{\sc A.~Makeev, D.~Maroudas, A.~Panagiotopoulos, and I.~Kevrekidis}, {\em
  Coarse bifurcation analysis of kinetic {Monte Carlo} simulations: a
  lattice-gas model with lateral interactions}, J. Chem. Phys., 117 (2002),
  pp.~8229--8240.

\bibitem{MatacheBS_NM00}
{\sc A.~Matache, I.~Babu$\breve{\textrm{s}}$ka, and C.~Schwab}, {\em
  Generalized {$p$-FEM} in homogenization}, Numer. Math., 86 (2000),
  pp.~319--375.

\bibitem{Milton02}
{\sc G.~Milton}, {\em Theory of composites}, Cambridge University Press,
  Cambridge, UK, 2002.

\bibitem{OberaiP_CMAME98}
{\sc A.~Oberai and P.~Pinsky}, {\em A multiscale finite element method for the
  {Helmholtz} equation}, Comput. Meth. Appl. Math. Engrg., 154 (1998),
  pp.~281--297.

\bibitem{RenardM_AWR97}
{\sc P.~Renard and G.~{De Marsily}}, {\em Calculating equivalent permeability:
  a review}, Adv. Water Res., 20 (1997), pp.~253--278.

\bibitem{MartinezGK_JCP04}
{\sc R.~{Rico-Martinez}, C.~Gear, and I.~Kevrekidis}, {\em Coarse projective
  {kMC} integration: forward/reverse initial and boundary value problems}, J.
  Comput. Phys., 196 (2004), pp.~474--489.
\newblock (original version can be found as nlin.CG/0307016 at arXiv.org).

\bibitem{RunborgTK_NL02}
{\sc O.~Runborg, C.~Theodoropoulos, and I.~Kevrekidis}, {\em Effective
  bifurcation analysis: a time-stepper based approach}, Nonlinearity, 15
  (2002), pp.~491--511.

\bibitem{SamaeyRK_SIMMS04}
{\sc G.~Samaey, D.~Roose, and I.~Kevrekidis}, {\em The gap-tooth scheme for
  homogenization problems}, SIAM J. Multiscale Model. Simulation, accepted
  (2004).
\newblock original version can be found as physics/0312004 at arXiv.org.

\bibitem{SirisupXKK_JCP05}
{\sc S.~Sirisup, D.~Xiu, G.~Karniadakis, and I.~Kevrekidis}, {\em
  Equation-free/{Galerkin}-free {POD}-assisted computation of incompressible
  flows}, J. Comput. Phys., in press (2005).

\bibitem{TheodoropoulosQK_PNAS00}
{\sc C.~Theodoropoulos, Y.~Qian, and I.~Kevrekidis}, {\em Coarse stability and
  bifurcation analysis using time-steppers: a reaction-diffusion example},
  Proc. Natl. Acad. Sci., 97 (2000), pp.~9840--9843.

\bibitem{XiuK_CMAME02}
{\sc D.~Xiu and G.~Karniadakis}, {\em Modeling uncertainty in steady state
  diffusion problems via generalized polynomial chaos}, Comput. Methods Appl.
  Math. Engrg., 191 (2002), pp.~4927--4948.

\bibitem{XiuK_SISC02}
\leavevmode\vrule height 2pt depth -1.6pt width 23pt, {\em The {Wiener-Askey}
  polynomial chaos for stochastic differential equations}, SIAM J. Sci.
  Comput., 24 (2002), pp.~619--644.

\bibitem{XiuK_JCP03}
\leavevmode\vrule height 2pt depth -1.6pt width 23pt, {\em Modeling uncertainty
  in flow simulations via generalized polynomial chaos}, J. Comput. Phys., 187
  (2003), pp.~137--167.

\bibitem{XiuK_IJHMT03}
\leavevmode\vrule height 2pt depth -1.6pt width 23pt, {\em A new stochastic
  approach to transient heat conduction modeling with uncertainty}, Inter. J.
  Heat Mass Trans., 46 (2003), pp.~4681--4693.

\end{thebibliography}

\end{document}